\documentclass[12pt]{amsart}
\usepackage[a4paper, left=3cm, right=3cm, height=22cm]{geometry}
\usepackage{graphicx}
\usepackage{amsmath,amsfonts,amsthm, amssymb}
\usepackage[arrow,matrix, all,cmtip]{xy}
\usepackage{enumerate}
\usepackage{cite}
\usepackage[usenames]{color}
\usepackage{amssymb}

\newtheorem{theorem}{Theorem}[section]
\numberwithin{equation}{section}
\newtheorem{prop}[theorem]{Proposition}

\newtheorem{cor}[theorem]{Corollary}

\theoremstyle{remark}
\newtheorem{df}[equation]{Definition}

\newtheorem{rmk}[theorem]{Remark}

\newcommand{\mb}{\mathbb}
\newcommand{\mf}{\mathfrak}

\newcommand{\mr}{\mathrm}

\newcommand{\C}{C} 

\usepackage[pdfauthor={Jin Cao}, %
				pdftitle={Cdga}%
			dvips,colorlinks=true]{hyperref}
\hypersetup{
	bookmarksnumbered=true,
	linkcolor=black,
}
\newcounter{elno}

\def\Z{\mathbb{Z}}                   
\def\Q{\mathbb{Q}}                   
\def\C{\mathbb{C}}                   
\def\uhp{{\mathbb H}}                
\def\dR{{\rm dR}}                    
\def\ker{{\rm ker}}              
\def\ker{{\rm ker}}              

\newcommand{\mat}[4]{
     \begin{pmatrix}
            #1 & #2 \\
            #3 & #4
       \end{pmatrix}
    }

\def\tmap{{\sf t}}

\begin{document}
\author{Jin Cao, Hossein Movasati, Roberto Villaflor Loyola}
\title{Gauss-Manin connection in disguise: Quasi Jacobi forms of index zero}
\date{\today}
\maketitle

\begin{abstract}
 We consider the moduli space of abelian varieties with two marked points and a frame of the relative de Rham cohomolgy with boundary at these points compatible with its mixed Hodge structure. Such a moduli space gives a natural algebro-geometric framework for higher genus quasi Jacobi forms of index zero and their differential equations which are given as vector fields. In the case of elliptic curves we compute explicitly the Gauss-Manin connection and such vector fields. 
\end{abstract}

\section{Introduction}
The literature on modular forms is a vast one, and has played a central role in number theory since its origins in the first half of the 19th century \cite{jacobi1829fundamenta}. The greatest achievement of the theory of modular forms has been the arithmetic modularity theorem \cite{taylor1995ring, breuil2001modularity} and the celebrated proof of Fermat's last theorem as a consequence of it \cite{wiles1995modular}. For an account of this part of the theory of modular forms, we refer the reader to the book \cite{diamond2005first}. In the modern formulation of the theory, one of the first steps is to interpret modular and automorphic forms as sections of line bundles over the so called modular curves, which are completions of quotients of the upper half plane $\mf{h}$ by the action of a congruence subgroup $\Gamma<\text{SL}_2(\Z)$. If one intends to generalize the classical theory of modular forms, it is natural to look for higher dimensional analogues of modular curves, these are the so called Shimura varieties. Roughly speaking, one of the big achievements of Hodge theory is to associate to every moduli space of algebraic varieties of a given type, a natural analytic variety $\Gamma\backslash\mathsf{D}$ given as the quotient of a homogeneous space $\mathsf{D}$, the so called Griffiths period domain, by the action of a discrete group $\Gamma$. This variety corresponds to a Shimura variety when $\mathsf{D}$ is a Hermitian symmetric domain. This is  the case only for few examples, for instance when  the Hodge structure of the underlying algebraic varieties is of weight 1 (curves and abelian varieties) or of weight 2 with $h^{2,0}=1$ (K3 surfaces, and other varieties with Hodge structure of level two like cubic fourfolds). In this classical setting, the theory of modular forms has found a fertile ground for generalizations, gaining a geometric framework for the theory of Siegel and Hilbert modular forms, and many types of automorphic forms on Hermitian symmetric domains. One of the reasons why this was possible, is due to the development of good compactifications of such spaces by Satake-Baily-Borel and the subsequent toroidal compactifications by Mumford. For a survey on these topics and some applications of the interplay between automorphic forms and moduli problems see \cite{laza2016perspectives}. For an introduction to the theory of Siegel and Hilbert modular forms and some of their applications to number theory see \cite{zagier123}. Recently, promising improvements has been obtained in the development of compactifications of the Griffiths period domain for non Hermitian symmetric cases, with applications to moduli problems. For a survey on these we refer to \cite{griffiths2018moduli}. Nevertheless, applications to the theory of modular and automorphic forms still seem to be far reaching for the non Hermitian symmetric case. 

The theory of modular forms can take another interesting direction if one looks for the differential equations which are satisfied by modular forms. This leads naturally to the concept of quasi modular forms. In fact, the algebra of quasi modular forms can be regarded as the smallest algebra closed under derivations which contains the algebra of modular forms. If one wants to develop a geometric framework for quasi modular forms and their differential equations, one realizes that the approach using Griffiths period domain is not the suitable one (even for the classical Hermitian symmetric domains). One of the main goals  of the project Gauss-Manin connection in disguise is to give an algebro-geometric framework for differential equations of quasi modular and automorphic forms, suitable for generalizations. The project started in \cite{M2012}, where the second author developed the geometric framework for quasi modular forms for $\text{SL}_2(\Z)$ and their differential equations, the so called Ramanujan equations. The main novelty in that work, was the introduction of a generalized period domain, having the classical Griffiths period domain as a quotient by the actions of an algebraic group. Since then, the project has been developed in many articles, providing the framework for several generalizations of modular forms and their differential equations. In \cite{ho18, Fonseca} were studied the quasi modular functions attached to moduli of abelian varieties, corresponding to Siegel and Hilbert quasi modular forms. One of the main achievements of the project Gauss-Manin connection in disguise is that it has been suitable to study moduli of Calabi-Yau varieties with Hodge structures of higher weights, where the period domain is not Hermitian symmetric. In \cite{ho22} the program is applied to the family of mirror quintic threefolds, obtaining a modular interpretation of the so called Yukawa coupling introduced by the physicists \cite{Candelas:1990rm} in the context of mirror symmetry. Later it was applied to the full family of Calabi-Yau threefolds \cite{HosseinMurad} resulting into a geometrization of topological string partition functions. Several other families of Calabi-Yau varieties have been studied providing new interesting functions with modular properties (see for instance \cite{MNDwork} for generalized Yukawa couplings attached to the Dwork family). In particular, the book \cite{M2020} contains the summary and the latest status of the project. This book was written in the most general framework of Hodge structures, and it is quite natural to rewrite it via the theory for mixed Hodge structures. In this paper we push forward the first step in this direction. 

Our aim is to develop the geometrization of another classical generalization of quasi modular forms, the so called quasi Jacobi forms of index zero. Jacobi forms can be thought as a cross between modular forms and elliptic functions. In fact, a Jacobi form of weight $k$ and index $m$ is a two complex variables function $\phi: \mb{H}_1 \times \mb{C} \to \mb{C}$ satisfying some functional equations relative to the discrete group $\text{SL}_2(\Z)\ltimes\Z$ and some Fourier expansion conditions (see Subsection \ref{sec9.1} for the index zero case). They were studied systematically first by Eichler and Zagier in \cite{EZ}. Their corresponding Hecke theory was developed, and their relations with other types of modular forms, such as Siegel modular forms, were investigated from an analytic viewpoint. 
Higher dimensional analogues of Jacobi forms appeared in several works \cite{shimura1978certain,gritsenko1984action, yamazaki1988jacobi, murase1991functions}, but the first attempt to build the general theory in the spirit of \cite{EZ} was done by Ziegler \cite{Ziegler1989}. On the other hand Kramer \cite{kramer1991geometrical} provided the first geometrization of higher level Jacobi forms as sections of a subsheaf of a line bundle over an elliptic modular surface $X_\Gamma$ associated to $\Gamma<\text{SL}_2(\Z)\ltimes\Z$ in the sense of \cite{shioda1972elliptic}. As in the case of quasi modular forms, quasi Jacobi forms arise when one looks for the smallest algebra closed under derivations containing the algebra of Jacobi forms. Recently, Libgober \cite{Libgober} was able to use the algebra of quasi Jacobi forms to study the elliptic genus of complex manifolds. This was done in analogy of the relation between quasi modular forms and the Witten genus, as pointed out by Zagier \cite{zagier1988note}. 
Besides the above, another source where quasi Jacobi forms have gained an increasing interest is a paper series by Oberdieck and Pixton \cite{Ober2018, Ober2019}, where they have found new connections between quasi Jacobi forms and Gromov-Witten invariants.


As mentioned before, in order to construct our geometrization of quasi Jacobi forms of index zero, we adapt the Gauss-Manin connection in disguise program to the framework of mixed Hodge theory. We start from the moduli space $\sf{T}$ of abelian varieties $X$ with two marked points $Y=\{O,P\}$  and a frame of the relative de Rham cohomology $H^1_\dR(X,Y)$ compatible with the mixed Hodge structure and the constant polarization (see Definition \ref{defT}). Our first result is the following:

\begin{theorem}
\label{theo1}
The moduli space $\mathsf{T}$ is a quasi-projective variety over $\Q$. In the case of elliptic curves, $\mathsf{T}$ is the affine variety
$$
\mathsf{T}=\mr{Spec} \ \C[a,b,c,t_1,t_2,\frac{1}{\Delta}]
$$
where $\Delta=27t_3^2-t_2^3$ and $t_3=4a^3-t_2a-b^2$. Moreover, $\mathsf{T}$ admits the universal family given by
$$
X=\{y^2=4x^3-t_2x-t_3\} \ , \ \ \ Y=\{O,P\} \ , \ O=(0:1:0), \ P=(a:b:1),
$$
and the frame of differential forms
$$
\mr{d}\left(\frac{x-a}{x}\right) \ , \ \  \frac{\mr{d}x}{y} \ , \ \   \left(c-\frac{b}{2a}\right)\mr{d}\left(\frac{x-a}{x}\right)+t_1\frac{\mr{d}x}{y} + \frac{x\mr{d}x}{y}-\mr{d}\left(\frac{y}{2x}\right).
$$
\end{theorem}

The algebra of regular functions in $\sf{T}$ are interpreted as quasi Jacobi forms of index zero. The differential equations of such quasi Jacobi forms are realized as vector fields (which we call modular vector fields) in $\sf{T}$ which in turn can be computed from the Gauss-Manin connection of the family of abelian varieties over $\sf{T}$. In terms of foliations, what one looks for are vector fields such that there exists some leaf parametrized by the generators of the algebra of quasi Jacobi forms of index zero. We proceed the other way around, by starting from a transcendental map $\tmap: \mathsf{D}\rightarrow\mathsf{T}$ defined Hodge theoretically (see \eqref{14ju2021} in Subsection \ref{sec9.1}) from the classical Griffiths period domain $\mathsf{D}$ to $\mathsf{T}$, and then we look for vector fields having the image of $\tmap$ as a leaf. In our case, $\mathsf{D}$ is naturally identified with an open set of $\mb{H}_g\times \C^g$ (see \eqref{Griffpd} in Subsection \ref{sec8.2}), and the modular vector fields on $\mathsf{T}$ are the algebraic incarnation of the derivations/vector fields $\frac{\partial}{\partial \tau_{ij}}, 1 \leq i \leq j \leq g$ and $\frac{\partial}{\partial z_k}, 1 \leq k \leq g$
in $(\tau,z)\in \mb{H}_g\times \C^g$. 

\begin{theorem}
\label{23june2021}
There are unique vector fields $\mathsf{v}_{ij}, 1 \leq i \leq j \leq g$ and $\mathsf{v}_k, 1 \leq k \leq g$ defined over $\mathbb{Q}$ in the moduli space $\mathsf{T}$ such that
\begin{equation}
A_{\mathsf{v}_{ij}} = C_{ij}
\end{equation}
and 
\begin{equation}
A_{\mathsf{v}_{k}} = C_{k},
\end{equation}
where $A_{\mathsf{v}_{ij}}$ (resp. $A_{\mathsf{v}_{k}}$) is the Gauss-Manin connection matrix composed with the vector field $\mathsf{v}_{ij}$ (resp. $\mathsf{v}_k$) and $C_{ij}$ (resp. $C_k$) is the constant matrix defined as above. Moreover, the Lie bracket of two such vector fields is zero.
\end{theorem}

In spite that the image of the map $\tmap$ is transcendental, the modular vector fields turn out being algebraic over $\mathsf{T}$. In the case of elliptic curves, using the explicit description of $\mathsf{T}$ given in Theorem \ref{theo1}, we compute explicitly the two modular vector fields $\mathsf{R}_z$ and $\mathsf{R}_\tau$ in terms of the parameters of $\mathsf{T}$. Moreover we show that:   

\begin{theorem} \label{23.06.2021}
There are unique global vector fields $\mathsf{R}_{\tau}$ and $\mathsf{R}_z$ on  $\mathsf{T}$ such that
\begin{equation} \label{R1=tau}
\nabla_{\mathsf{R}_{\tau}} \left( \begin{matrix}
\alpha_1 \\
\alpha_2 \\
\alpha_3
\end{matrix}\right) = \left( \begin{matrix}
0 & 0 & 0 \\
0 & 0 & -1 \\
0 & 0 & 0
\end{matrix}\right)\left( \begin{matrix}
\alpha_1 \\
\alpha_2 \\
\alpha_3
\end{matrix}\right) 
\end{equation}
and 
\begin{equation} \label{R2=z}
\nabla_{\mathsf{R}_z} \left( \begin{matrix}
\alpha_1 \\
\alpha_2 \\
\alpha_3
\end{matrix}\right) = \left( \begin{matrix}
0 & 0 & 0 \\
-1 & 0 & 0 \\
0 & 0 & 0
\end{matrix}\right)\left( \begin{matrix}
\alpha_1 \\
\alpha_2 \\
\alpha_3
\end{matrix}\right) ,
\end{equation}
where $\nabla$ is the Gauss-Manin connection of $\mathsf{T}$. More precisely, we let $t_3 = 4a^3 - t_2a -b^2$ and then
\begin{equation}
\begin{split}
\mathsf{R}_{\tau} &= (-2a^2+2at_1+bc+\frac{t_2}{3})\frac{\partial}{\partial a}+(6a^2c-\frac{ct_2}{2}-3ab+3bt_1)\frac{\partial}{\partial b} \\
& +(ac+ct_1-\frac{b}{2})\frac{\partial}{\partial c} +(t_1^2-\frac{t_2}{12})\frac{\partial}{\partial t_1} + (4t_1t_2-6t_3)\frac{\partial}{\partial t_2} 
\end{split}
\end{equation}
\begin{equation}
\mathsf{R}_z = b\frac{\partial}{\partial a}+(6a^2-\frac{t_2}{2})\frac{\partial}{\partial b}+(a+t_1)\frac{\partial}{\partial c}.
\end{equation}
\end{theorem}

The translation into holomorphic context of Jacobi forms is done through the $\tmap$ map \eqref{14ju2021} defined in Subsection \ref{sec9.1}. In this way the coordinate functions $a,b,c,t_1,t_2$ are transformed up to constants into Jacobi forms 
$\wp(\tau, z), \wp^{'}(\tau, z), J_1(\tau, z),E_2(\tau), E_4(\tau)$ respectively. More precisely we prove the following:

\begin{theorem}
\label{2021china}
The pullback of $a,b, c,t_1,t_2,t_3$ under the map $\tmap$ in \eqref{14ju2021}  
are 
$$
(2\pi i)\wp(\tau, z) \ , \ \  (2\pi i)^{\frac{3}{2}}\wp^{'}(\tau, z) \ , \ \ - (2\pi i)^{\frac{1}{2}}J_1(\tau, z),
$$
$$
-\frac{2\pi i}{12}E_2(\tau)\ , \ \ 12\left(\frac{2\pi i}{12}\right)^2E_4(\tau) \ , \ \ -8\left(\frac{2\pi i}{12}\right)^3E_6(\tau),
$$
respectively.
\end{theorem}

The vector fields $\mathsf{R}_\tau$ and $\mathsf{R}_z$ as differential equations between these Jacobi forms correspond to those already computed in \cite[Lemma 48]{Ober2018}.

Since we want to emphasize on the algebraicity of the modular vector fields, we devote the first part of the article to the purely algebraic results and constructions. Only the last two sections are devoted to the transcendental constructions and their relation with the algebraic ones. Keeping this in mind the text is organized as follows. In Section \ref{sect2} we recall the algebraic definition of the relative de Rham cohomology of a pair. Section \ref{sect3} is devoted to the computation of the cup product in relative de Rham cohomology. In Section \ref{section4} we describe the polarized mixed Hodge structure of the relative de Rham cohomology of an abelian variety relative to two points. Using this description we give the precise definition of the moduli space $\mathsf{T}$ in Section \ref{sect5}, and we prove Theorem \ref{theo1}. In this section we also describe a natural action of an algebraic group $\mathsf{G}$ on $\mathsf{T}$. This action will be key to define the $\tmap$ map, and to describe the modular functional equations satisfied by the solutions of the modular vector fields after lifting them via $\tmap$. Section \ref{section6} is devoted to the computation of the Gauss-Manin connection on the relative de Rham cohomology bundle of the universal family of $\mathsf{T}$. Using the Gauss-Manin connection computations, we prove Theorem \ref{23.06.2021} in Section \ref{sect7}. In Section \ref{gpd2021} we shift to the transcendental objects. In order to define later the map $\tmap$ in terms of the Hodge theoretic information of $\mathsf{T}$, we start by defining the generalized period domain $\mathsf{\Pi}$ in our context. As for the classical Griffiths period domain, it comes with a natural action of a discrete group $\Gamma_\Z$. Moreover, the algebraic group $\mathsf{G}$ described in Section \ref{sect5} also acts on $\mathsf{\Pi}$, in such a way that $\mathsf{\Pi}/\mathsf{G}$ corresponds to the classical Griffiths period domain $\mathsf{D}$. We introduce the period map $\mathsf{P}:\mathsf{T}\rightarrow \mathsf{U}=\Gamma_\Z\backslash \mathsf{\Pi}$ and show it is a biholomorphism (see Proposition \ref{isoperiodmap}). We define also the $\tau$-map as a section of the quotient map $\mathsf{\Pi}\rightarrow\mathsf{D}$ (see Subsection \ref{23j2021}). Using all these constructions we justify our choices of the constant matrices $C_{ij}$ and $C_k$ of Theorem \ref{23june2021}, and complete its proof. Finally in Section \ref{sect9} we define the map $\tmap$ as the composition of the $\tau$-map and the inverse of the period map. Using it we translate the differential equations given by the modular vector fields $\mathsf{R}_\tau$ and $\mathsf{R}_z$ into the classical differential equations of quasi Jacobi forms of index zero, thus completing the proof of Theorem \ref{2021china}. Throughout the text $k$ is a field of characteristic $0$. In some case we will consider  it as a subfield of the field of complex numbers $\mb{C}$.



\section{Relative de Rham cohomology}
\label{sect2}

For the sake of completeness, in this section we recall the algebraic definition of the relative de Rham cohomology in terms of hypercohomology. In the case of elliptic curves with two marked points we give an alternative description in terms of global meromorphic forms without residues. In that case we also provide a explicit basis for both descriptions. The hypercohomology description will be useful when we compute the cup product in the next section. The 
meromorphic description will be used later for the computation of the Gauss-Manin connection (see Section \ref{section6}).
\subsection{The definition}
\begin{df} \label{relative cohomology}
Let $X$ be a smooth variety over $k$ and $Y$ be a smooth subvariety of $X$. We consider the complex $(\Omega^{\bullet}_{X/k}, d)$ (resp. $(\Omega^{\bullet}_{Y/k}, d)$) of regular differential forms on $X$ (resp. $Y$).
The (algebraic) relative de Rham cohomology of $(X, Y)$ is defined to be the hypercohomology of the following complex
\[
H^m_{\mr{dR}}((X, Y)/k) := \mb{H}^m(\Omega_{(X, Y)/k}^{\bullet}, d),
\]
where
\[
\Omega^m_{(X, Y)/k} := \Omega^m_{X/k} \oplus \Omega_{Y/k}^{m-1}
\]
and
\[
d: \Omega^m_{(X, Y)/k} \to \Omega^{m+1}_{(X, Y)/k} \ , \  (\omega, \alpha) \to (d\omega, \omega|_Y - d\alpha).
\]
\end{df}
The $C^{\infty}$ relative de Rham cohomology of a pair of manifolds can be defined similarly, see \cite[page 78]{BottTu}. In fact, given two manifolds $X^{\infty}, Y^{\infty}$ together with a continuous map between manifolds $f: Y^{\infty} \to X^{\infty}$, one can define the complex
\[
\Omega^{\bullet}(f) := \Omega^{\bullet}(X^{\infty}) \oplus \Omega^{\bullet - 1}(Y^{\infty}),
\]
where $\Omega^{\bullet}$ is the complex of $C^{\infty}$ differential forms,
with differential 
\[
d(\omega, \alpha) _= (d\omega, f^*\omega - d\alpha).
\]
The relative de Rham cohomology for $(X^{\infty}, Y^{\infty}, f: Y^{\infty} \to X^{\infty})$ is defined to be the cohomology of $\Omega^{\bullet}(f)$:
\[
H^m_{\mr{dR}}(X^{\infty}, Y^{\infty}) := H^m(\Omega^{\bullet}(f), d).
\]
In particular, if $f$ is an immersion $i$ of a submanifold $Y^{\infty}$ inside $X^{\infty}$, we write $\Omega^{\bullet}(i)$ simply as $\Omega^{\bullet}_{(X^{\infty}, Y^{\infty})}$.

\begin{prop}
Assume that $k = \mb{C}$. Let $X$ be a smooth variety over $k$ and $Y$ be a smooth subvariety of $X$. For any $m \ge 0$, we have a natural isomorphism
\[
H^m_{\mr{dR}}((X, Y)/\mb{C}) \cong H^m_{\mr{dR}}(X^{\infty}, Y^{\infty}).
\]
\end{prop}
\begin{proof}
Note that we have  the following commutative diagram of complexes
\begin{equation}
\begin{xy}
(60,0)*+{0}="v1";(90,0)*+{\Omega^{\bullet-1}_{Y/\mb{C}}}="v2";
(120,0)*+{\Omega^{\bullet}_{(X, Y)/\mb{C}}}="v3";
(150,0)*+{\Omega^{\bullet}_{X/\mb{C}}}="v4";(180,0)*+{0}="v5";
(60,-20)*+{0}="v6";(90,-20)*+{\Omega^{\bullet-1}_{Y^{\infty}}}="v7";
(120,-20)*+{\Omega^{\bullet}_{(X^{\infty}, Y^{\infty})}}="v8";
(150,-20)*+{\Omega^{\bullet}_{X^{\infty}}}="v9";(180,-20)*+{0}="v10";
{\ar@{->} "v1";"v2"};{\ar@{->} "v2";"v3"};{\ar@{->} "v3";"v4"};{\ar@{->} "v4";"v5"};
{\ar@{->} "v6";"v7"};{\ar@{->} "v7";"v8"};{\ar@{->} "v8";"v9"};{\ar@{->} "v9";"v10"};
{\ar@{->} "v2";"v7"};{\ar@{->} "v3";"v8"};{\ar@{->} "v4";"v9"};
\end{xy}
\end{equation}
where the vertical maps are just the natural inclusion of sheaves with respect to the Zariski topology. By Grothendieck's algebraic description of the de Rham cohomology (or equivalently by Atiyah-Hodge theorem, c.f. \cite[Section 5.3]{ML}) the left and right vertical map are quasi isomorphisms over the Zariski topology. Therefore the middle vertical map is also a quasi isomorphism. Since the hypercohomolgy of $\Omega_{(X^\infty, Y^\infty)}^\bullet$ is the same over the Zariski topology than over the analytic topology, the result follows noting that $\Omega_{(X^\infty, Y^\infty)}^\bullet$ is fine over the analytic topology (see for instance \cite[Proposition 3.2, Proposition 3.6]{ML}).
\end{proof}

\subsection{Relative de Rham cohomology for abelian varieties}
\label{10may2020}
Consider an abelian variety $X$ of dimension $g$, and $Y$ as a subvariety of $X$ consists of two points $O$ and $P$.   The short exact sequence
\[
0 \to  \Omega^{\bullet-1}_{Y/k} \to \Omega^{\bullet}_{(X,Y)/k} \to \Omega^{\bullet}_{X/k} \to 0
\]
induces the long exact sequence
\[
\cdots \to \mb{H}^0(X, \Omega^{\bullet}_{X/k}) \to \mb{H}^1(Y,  \Omega^{\bullet-1}_{Y/k} ) \to \mb{H}^1(X, \Omega^{\bullet}_{(X,Y)/k}) \to \mb{H}^1(X, \Omega^{\bullet}_{X/k}) \to \cdots
\]
which translates into the exact sequence
\begin{equation}
\label{11may2021}
0 \to H^0_{\mr{dR}}(X/k) \to H^0_{\mr{dR}}(Y/k) \to H^1_{\mr{dR}}((X, Y)/k) \to H^1_{\mr{dR}}(X/k) \to 0,
\end{equation}
together with the isomorphisms
\[
H^i_{\mr{dR}}((X, Y)/k) \cong H^i_{\mr{dR}}(X/k), \ \ \ \ \forall  i \geq 2.
\]
It follows from \eqref{11may2021} that $\dim H^1_\dR((X,Y)/k)=2g+1$. 
One can easily verify that  
$$
\text{coker }(H^0_\dR(X/k)\rightarrow
H^0_\dR(Y/k))= k\cdot f,
$$
where $f:Y\rightarrow k$ is given by $f(O)=1$ and $f(P)=0$.

\subsection{Relative de Rham cohomology for elliptic curves}
 \label{computation of the relative cohomology}
Using \v{C}ech complexes we can compute the relative cohomology. Let us consider the following case: Let $X$ be a smooth projective curve of genus one given by the Weierstrass equation 
\[
X=\{y^2 = 4x^3 - t_2x -t_3\}
\] 
and let $Y := \{O, P\}$, where $O$ is the infinity point and $P = (a, b)$. Hence $t_3=4a^3-t_2a-b^2$. Then we can choose the affine open cover 
\[
U_0 = \mr{Spec}(\mb{C}[x, y] / (y^2 - 4x^3 +t_2x +t_3)) = X-O = X - \{[0:1:0]\}
\]
and 
\[
U_1 = (X - \{x = 0\})\cup\{O\}.
\]
According to the definition, we compute the relative cohomology via the following complex
\begin{equation}
\begin{split}
&0 \to \Omega^0(U_0)\oplus \Omega^0(U_1) \to \Omega^0(U_0 \cap U_1) \oplus \Omega^1(U_0, U_0 \cap Y) \oplus \Omega^1(U_1, U_1\cap Y)   \\
&\to \Omega^1(U_0 \cap U_1, U_0 \cap U_1 \cap Y) \oplus \Omega^2(U_0) \oplus \Omega^2(U_1) \to \Omega^2(U_0 \cap U_1) \to 0.
\end{split}
\end{equation}
Here for the differential of the double complex, we choose the sign rule defined in \cite[page 29]{ML}, this is $D=\delta+(-1)^q\mr{d}$.
In particular, we have
\begin{equation}
\begin{split}
& H^1_{\mr{dR}}(X, Y) = \\
&\frac{\{((\omega_0, \alpha_0), (\omega_1, \alpha_1), f_{01} )\mid \mr{d}f_{01} = \omega_1 | _{U_{01}} - \omega_0 |_{U_{01}} \ , \ f_{01}|_{Y}=\alpha_1|_{U_{10}}-\alpha_0|_{U_{10}} \}}{\{( (\mr{d}f_0, f_0|_{Y\cap U_0}), (\mr{d}f_1, f_1|_{Y\cap U_1}), f_1|_{U_{01}}- f_0|_{U_{01}})\}}
\end{split}
\end{equation}
where $f_i \in \Omega^0(U_i), f_{01} \in \Omega^0(U_0 \cap U_1), (\omega_i, \alpha_i) \in \Omega_{(X,Y)}^1(U_i, U_i \cap Y)$ and $U_{01} = U_0 \cap U_1$.

\begin{prop} \label{basis of relative cohomology}
For $a \not=0$, we can choose a basis of $H^1_{\mr{dR}}((X, Y)/\C)$ as follows:
\begin{enumerate}
\item
$((0, 0), (\mr{d}f, 0), f|_{U_{01}})$, where $f= \frac{x-a}{x}$;
\item
$((\frac{\mr{d}x}{y}|_{U_0}, 0), (\frac{\mr{d}x}{y}|_{U_1}, 0), 0)$, where $\frac{\mr{d}x}{y}$ is a holomorphic $1$-form on $X$;
\item
$((\frac{x\mr{d}x}{y}, 0), (\frac{x\mr{d}x}{y}+\mr{d}g, g|_Y), g|_{U_{01}})$, where $g = -\frac{y}{2x}$.
\end{enumerate}
\end{prop}
\begin{proof}
For the second and third item, they form a basis of $H^1_{\mr{dR}}(X/\C)$ and the details can be found in \cite[Proposition 2.4]{M2012}. For the first one, it is enough to show that this element is not zero in $H^1_{\mr{dR}}((X, Y)/\C)$. If this is not true, then we can write $( (0, 0), (\mr{d}f, 0), f|_{U_{01}})$ as
\[
((\mr{d}f_0, f_0|_{Y\cap U_0}), (\mr{d}f_1, f_1|_{Y\cap U_1}), f_1|_{U_{01}}- f_0|_{U_{01}}).
\]
Then $f_0 = 0$ and hence $f_1=f$. However the infinite point $O \in Y \cap U_1$ and $f(O) = 1$, which is a contradiction with $f_1|_{Y\cap U_1}=0$.
\end{proof}

\begin{rmk}
For $a = 0$, we have $df=0$ and so the basis we have chosen above  degenerates. 
\end{rmk}

\subsection{Meromorphic forms without residues} In this section we provide another algebraic interpretation of the relative de Rham cohomology, valid only for the case of elliptic curves. It depends on the choice of an affine chart $U$ containing $Y=\{O,P\}$. The advantage of this description is that the Gauss-Manin connection becomes much simpler to compute (see Subsection \ref{subsectalgGMell}).

\begin{prop}
\label{eq2.6}
Let $U\subseteq X$ be an affine open set such that $Y\subseteq U$. Then 
$$
H^1_\dR(X,Y)\cong \frac{\text{$\omega\in \Gamma(\Omega^1_U)$ without residues on $X- U$}}{\text{exact forms $df$ with $f|_Y=0$}}.
$$
\end{prop}

\begin{proof}
Recall that
$$
H^1_\dR(X,Y)=\uhp^1(X,\Omega_{(X,Y)}^\bullet).
$$
Note first that the complex $\Omega_{(X,Y)}^\bullet$ is a resolution of the sheaf
$$
\mathcal{K}:=\ker(\C_X\xrightarrow{|_Y} \C_Y).
$$
Consider the following complex of sheaves over $X$
$$
\hat{\Omega}_U^\bullet: \hspace{0.5cm} i_*\mathcal{I}_Y\xrightarrow{d}i_*\Omega_U^1\xrightarrow{Res} \C_{X- U}\rightarrow 0 \ ,
$$
where $\mathcal{I}_Y$ is the ideal sheaf of $Y$ over $U$, and $i:U\hookrightarrow X$ is the inclusion map. It is also easy to see that $\hat{\Omega}_U^\bullet$ is also a resolution of $\mathcal{K}$. Therefore
$$
H^1_\dR(X,Y)=\uhp^1(X,\Omega_{(X,Y)}^\bullet)\cong H^1(X,\mathcal{K})\cong\uhp^1(X,\hat{\Omega}_U^\bullet).
$$
Finally, note that $\hat{\Omega}_U^\bullet$ is acyclic since each sheaf is supported on the affine set $U$ or in the finite set $X- U$, thus
$$
\uhp^1(X,\hat{\Omega}_U^\bullet)\cong H^1(\Gamma(\hat{\Omega}_U^\bullet)).
$$
\end{proof}

\begin{prop}
\label{isoeq2.6}
Under the same hypothesis of Proposition \ref{eq2.6}, let $Y=\{O,P\}$, $V= X-O$, $f\in \Gamma(\mathcal{O}_U)$ be such that $f(P)=0$ and $f(O)=1$, and take $\mathcal{U}=\{V,U\}$, then the map
\begin{equation}
\label{22june2021}
\uhp^1(\mathcal{U},\Omega_{(X,Y)}^\bullet)\rightarrow H^1(\Gamma(\hat{\Omega}_U^\bullet))
\end{equation}
$$
\omega=((\omega_0,f_0),(\omega_1,f_1),f_{01})\mapsto \omega_1+(f_1(P)-f_1(O))df
$$
is an isomorphism. This isomorphism is independent of the choice of $f$.
\end{prop}

\begin{proof}
First of all, to see that it is well-defined we have to show that $\omega_1$ has no residues. In fact, since $\omega_1-\omega_0=df_{01}$, it follows that $\omega_1$ has the same residues as $\omega_0$. On the other hand $\omega_0$ has at most one pole at $O$, hence it has no residues. Now, to show that it is an isomorphism, it is enough to show that it is injective. In fact, suppose that there exists some $h\in \mathcal{I}_Y(U)$ such that
$$
\omega_1+(f_1(P)-f_1(O))df=dh,
$$
then $\omega_1=d\eta$, for $\eta=(f_1(O)-f_1(P))f+h$. Replacing $\omega$ by
$$
\omega-D(0,\eta+f_1(P))=((\omega_0,f_0),(0,0),f_{01}-\eta-f_1(P))
$$
it follows that $\omega_0=d\mu$, for $\mu=-f_{01}+\eta+f_1(P)$ and so $\mu$ only has poles at $O$, i.e. it is defined on $V$, then
$$
\omega-D(\mu,\eta+f_1(P))=((0,f_0-\mu|_{Y\cap V}),(0,0),0).
$$
Finally since $P\in U\cap V$ and $Y\cap V=\{P\}$ it follows that $f_0-\mu|_{Y\cap V}=0$, i.e. $\omega=D(\mu,\eta+f_1(P))$.
\end{proof}

\begin{rmk}
Under the isomorphism \eqref{22june2021}, we may choose the basis of the relative de Rham cohomology directly rather than using the representatives in the \v{C}ech complex. This will be useful when we compute the Gauss-Manin connection later. Depending on the coordinates of the point $P$, we choose the following basis on $H^1(\Gamma(\hat{\Omega}_U^\bullet))$: 
\begin{equation}
\label{24june2021}
\mr{d}\left(\frac{x-a}{x}\right), \frac{\mr{d}x}{y}, \frac{x\mr{d}x}{y} - \mr{d}\left(\frac{y}{2x}\right)-\frac{b}{2a}\mr{d}\left(\frac{x-a}{x}\right),\ \ \  a\neq0,
\end{equation}
\[
\mr{d}\left(\frac{x-a}{x-1}\right), \frac{\mr{d}x}{y}, \frac{x\mr{d}x}{y} - \mr{d}\left(\frac{y}{2(x-1)}\right)-\frac{b}{2(a-1)}\mr{d}\left(\frac{x-a}{x-1}\right),\ \ \  a\neq 1.
\]
Note that for $a\neq 0$ we are considering $U=(X-\{x=0\})\cup\{O\}$, while for $a\neq 1$ we take $U=(X-\{x=1\})\cup\{O\}$.   The first differential form $\omega_1=df$ is chosen in such a way that $f(O)=1, f(P)=0$. Note also that the correction of $\frac{xdx}{y}$ with an exact differential form kills its pole at $O$. The computations are similar as in \cite[Section 2.8]{M2012}. We remark that when $a\neq 0,1$ both  basis are equal in $H^1(\Gamma(\hat{\Omega}_U^\bullet))$ for $U=(X-\{x(x-1)=0\})\cup\{O\}$. For instance, the difference of the third element in both basis is an exact differential form $dg$ with $g(P)=g(O)=\frac{b}{2a}-\frac{b}{2(a-1)}$.  
\end{rmk}

\section{Relative cup product}
\label{sect3}

In this section we compute the cup product of forms in the algebraic relative de Rham cohomology described as elements of the hypercohomology of the complex of relative forms (see Definition \ref{relative cohomology}). This is necessary in order to compute the cup product of the basis given in Proposition \ref{basis of relative cohomology}. 
\subsection{Relative cup product in hypercohomology} At first we define the cup product in the usual relative de Rham cohomology. The basic idea comes from \cite[Section 3]{EV}.
\begin{df}
We define the product
\[
\cup: \Omega^m_{(X^{\infty}, Y^{\infty})} \otimes \Omega^n_{(X^{\infty}, Y^{\infty})}  \to \Omega^{m+n}_{(X^{\infty}, Y^{\infty})} 
\]
by
\[
(\check{\omega}, \check{\alpha}), (\check{\nu}, \check{\beta}) \to (\check{\omega} \wedge \check{\nu}, (-1)^m \check{\omega} |_{Y^{\infty}} \wedge  \check{\beta}).
\]
\end{df}
\begin{prop}
For $(\check{\omega}, \check{\alpha}) \in \Omega^m_{(X^{\infty}, Y^{\infty})}$ and $ (\check{\nu}, \check{\beta}) \in \Omega^n_{(X^{\infty}, Y^{\infty})}$, we have
\begin{equation} \label{boudary formula for cup product}
d((\check{\omega}, \check{\alpha}) \cup (\check{\nu}, \check{\beta})) = d((\check{\omega}, \check{\alpha}) ) \cup (\check{\nu}, \check{\beta}) + (-1)^m  (\check{\omega}, \check{\alpha})   \cup d((\check{\nu}, \check{\beta}) ).
\end{equation}
\end{prop}
\begin{proof}
This is just a straightforward computation
\begin{equation} \nonumber
\begin{split}
& \mr{LHS} = d((\check{\omega} \wedge \check{\nu}, (-1)^m \check{\omega} |_{Y^{\infty}} \wedge  \check{\beta})) = (d(\check{\omega} \wedge \check{\nu}), (\check{\omega} \wedge \check{\nu})|_{Y^{\infty}}- d((-1)^m \check{\omega} |_{Y^{\infty}} \wedge  \check{\beta})) = \\
& (d\check{\omega} \wedge \check{\nu} + (-1)^m \check{\omega} \wedge d \check{\nu}, (\check{\omega} \wedge \check{\nu})|_{Y^{\infty}} + (-1)^{m+1}d(\check{\omega} |_{Y^{\infty}}) \wedge \check{\beta} + (-1)^{2m+1} \check{\omega} |_{Y^{\infty}} \wedge d\check{\beta})
\end{split} 
\end{equation}
On the other hand, 
\begin{equation} \nonumber
\begin{split}
& \mr{RHS} = (d\check{\omega}, \check{\omega}|_{Y^{\infty}} - d\check{\alpha}) \cup (\check{\nu}, \check{\beta}) + (-1)^m (\check{\omega}, \check{\alpha})   \cup (d\check{\nu}, \check{\nu}|_{Y^{\infty}} - d\check{\beta})=  \\
& (d\check{\omega} \wedge \check{\nu}, (-1)^{m+1} d\check{\omega}|_{Y^{\infty}} \wedge \check{\beta}) + (-1)^m (\check{\omega} \wedge d\check{\nu}, (-1)^m \check{\omega}|_{Y^{\infty}} \wedge (\check{\nu}|_{Y^{\infty}} - d\check{\beta} )) = \\
& (d\check{\omega} \wedge \check{\nu} + (-1)^m \check{\omega} \wedge d \check{\nu}, (\check{\omega} \wedge \check{\nu})|_{Y^{\infty}} + (-1)^{m+1}d(\check{\omega} |_{Y^{\infty}}) \wedge \check{\beta} + (-1)^{2m+1} \check{\omega} |_{Y^{\infty}} \wedge d\check{\beta}).
\end{split} 
\end{equation}
\end{proof}
\begin{cor}
The product $\cup$ induces a product structure on the relative de Rham cohomology
\[
\cup: H^{m}_{\mr{dR}}(X^{\infty}, Y^{\infty}) \otimes H^{n}_{\mr{dR}}(X^{\infty}, Y^{\infty}) \to H^{m+n}_{\mr{dR}}(X^{\infty}, Y^{\infty}),
\]
which is called the cup product for the relative de Rham cohomology.
\end{cor}
Next we want to find the corresponding bilinear map in algebraic relative de Rham cohomology. According to Leray's theorem
, we may take a cover consisting of open affine subsets of $X$ to compute the algebraic relative de Rham cohomology in terms of a \v{C}ech resolution (see for instance \cite[Theorem 3.1]{ML}). Fix an affine open cover $\{U_i\}$ of $X$. For an element $(\omega, \alpha) \in H^m_{\mr{dR}}((X, Y)/k)$, it can be represented as the sum of $(\omega^{m-r}_{i_0i_1\cdots i_r}, \alpha^{m-1-r}_{i_0i_1\cdots i_r})$ for $r = 0, \ldots, m$. Here the upper index $m-r$ denotes $\omega^{m-r}_{i_0i_1\cdots i_r}$ is a $(m-r)$-forms and the lower index $i_0i_1\cdots i_r$ means the differential form lies in $U_{i_0} \cap U_{i_1} \cap \cdots\cap U_{i_r}$. Note that for $\alpha$, the lower index denotes the affine cover restricted to $Y$.
\begin{prop}
\label{propcuphyp}
Keeping the same notation as above, the cup product of $(\omega, \alpha) \in H^m_{\mr{dR}}((X, Y)/k)$ and $(\nu, \beta) \in H^n_{\mr{dR}}((X, Y)/k)$ is given by $(\mu, \gamma)$, where
\[
\mu^{n+m-j}_{i_0i_1\cdots i_j} = \sum^j_{r=0} (-1)^{m(j-r)+r(j-1)} \omega^{m-r}_{i_0\cdots i_r} \wedge \nu^{n-j+r}_{i_r\cdots i_j}
\]
for $0 \leq j \leq n+m$ and
\[
\gamma^{n+m-k-1}_{i_0i_1\cdots i_k} = (-1)^m\sum^k_{s=0} (-1)^{m(k-s)+s(k-1)} {\omega^{m-s}_{i_0\cdots i_s}}|_Y \wedge \beta^{n-k-1+s}_{i_s\cdots i_k}
\]
for $0 \leq k \leq n+m-1$.
\end{prop}
\begin{proof}
The proof is just the definition together with the usual twisting cup product formula in hypercohomology (see for instance \cite[Theorem 5.3]{ML}).
\end{proof}

\subsection{The case of elliptic curves}
Let us apply Proposition \ref{propcuphyp} to compute the cup of the elements of the basis given in Proposition \ref{basis of relative cohomology}. We choose the affine open cover $\{U_0, U_1, U_{01}\}$ of the smooth projective curve as in Subsection \ref{computation of the relative cohomology}. Take two elements $(\omega, \alpha)$ and $(\nu, \beta)$ which are represented as
\[
((\omega_0, \alpha_0), (\omega_1, \alpha_1), \omega_{01}), ( (\nu_0, \beta_0), (\nu_1, \beta_1), \nu_{01}).
\]
Then using the above formula, we can compute $(\omega, \alpha) \cup (\nu, \beta)$, which can be represented as
\[
(\omega_0 \wedge \nu_0, -\omega_0|_Y \wedge \beta_0), (\omega_1 \wedge \nu_1,  -\omega_1|_Y \wedge \beta_1), (-\omega_0 \wedge \nu_{01}+\omega_{01} \wedge \nu_1, -\omega_{01}|_Y \wedge \beta_1)
\]
In particular, using the basis $\omega_i, i = 1,2,3$ of $H^1_{\mr{dR}}((X, Y)/\C)$ given in Proposition \ref{basis of relative cohomology}, we get that
\begin{equation}
\omega_1 \cup \omega_2 = \left((0, 0), \left( \mr{d}f \wedge \frac{\mr{d}x}{y}, 0 \right), \left(\frac{(x-a)\mr{d}x}{xy}, 0\right)\right), 
\end{equation}
\begin{equation}
\omega_1 \cup \omega_3 = \left((0, 0), \left( \mr{d}f \wedge \frac{x\mr{d}x}{y}, 0\right), \left(\frac{(x-a)\mr{d}x}{y}, 0\right)\right), 
\end{equation}
\begin{equation}
\omega_2 \cup \omega_3 =\left( (0, 0), (0, 0), \left(-\frac{\mr{d}x}{2x}, 0\right) \right)
\end{equation}

\begin{rmk}
When $X$ is a smooth projective curve and $Y$ is a set of distinct points, we have the following commutative diagram
\begin{equation}
\label{tracediagram}
\begin{xy}
(60,0)*+{H^{1}_\dR(X, Y) \times H^{1}_\dR(X, Y)}="v1";(120,0)*+{H^{2}_\dR(X, Y)}="v2";
(60,-20)*+{H^{1}_\dR(X) \times H^{1}_\dR(X)}="v3";(120,-20)*+{H^{2}_\dR(X)}="v4";(160,-20)*+{\mb{C}}="v5";
{\ar@{->}^{\cup} "v1";"v2"};{\ar@{->}^{\cup} "v3";"v4"};{\ar@{->}^{\mathrm{Tr}} "v4";"v5"};
{\ar@{->}^{j^* \times j^*} "v1";"v3"};{\ar@{->}^{j^*} "v2";"v4"}
\end{xy}
\end{equation}
Here $j: (X, \emptyset) \to (X, Y)$ is the inclusion and note that both the right vertical map and the trace map are isomorphisms. Hence the composition of the relative cup product, $j^*$ and the trace map gives us the bilinear map
\[
\langle \cdot, \cdot \rangle: H^{1}_\dR((X, Y)/\C) \times H^{1}_\dR((X, Y)/\C) \to \mb{C},
\]
which is called the relative trace map. We apply the above constructions to our case, i.e. $X$ a smooth projective curve of genus one and $Y$ a subset of two distinct points. Composing the relative trace map with the relative cup product, we find that
\[
\langle \omega_2, \omega_3 \rangle =  - \langle \omega_3, \omega_2 \rangle = 1, 
\]
and the others are zero. Here we used the fact that the trace map is the residue of $\omega_{01}$ around the infinite point (see for instance \cite[Page 19]{M2012}).
\end{rmk}

\section{Mixed Hodge structure on relative cohomology}
\label{section4}
In this section, we briefly recall the construction of the mixed Hodge structure on the relative de Rham cohomology. For more details we refer to \cite[Definition 3.13, Theorem 3.18]{PS}.

\subsection{Polarized mixed Hodge structure for abelian varieties} As shown in \cite[Example 3.24]{PS}, one may use explicit representatives to define the mixed Hodge structure on relative cohomology. 
We want to explore this in detail in our case, that is, 
 $X$ is an abelian variety of dimension $g$ over $k$ and $Y$ is a closed subvariety of $X$ consisting of two distinct points. Using the exact sequence \eqref{11may2021} we can determine the mixed Hodge structure of
$H^1_{\mr{dR}}((X, Y)/k) $. In fact, the weight filtration of $H^1_{\mr{dR}}((X, Y)/k)$ corresponds to
\[
W_0 H^1_{\mr{dR}}((X, Y)/k) = \mr{coker}(H^0_{\mr{dR}}(X/k) \to H^0_{\mr{dR}}(Y/k))
\]
\[
\subset
W_1 H^1_{\mr{dR}}((X, Y)/k) = H^1_{\mr{dR}}((X, Y)/k),
\]
and the only nontrivial piece  $F^1$ of the Hodge filtration is given by
\begin{equation}
\begin{split}
0 \to &  F^1 (\mr{coker}(H^0_{\mr{dR}}(X/k) \to H^0_{\mr{dR}}(Y/k)))=0 \\
\to & F^1 H^1_{\mr{dR}}((X, Y)/k) \to  F^1 H^1_{\mr{dR}}(X/k) \to 0
\end{split}
\end{equation}
The following proposition computes a basis compatible with the polarized mixed Hodge structure of a polarized abelian variety with two marked points. This is a slight modification of \cite[Proposition 11.1]{M2020} into our context.
\begin{prop}
\label{prop1}
We can take a basis $\alpha_0,\alpha_1,\ldots,\alpha_{2g}\in H^1_\dR((X,Y)/k)$ such that $\alpha_0\in H^1_\dR((X,Y)/k)$, $\alpha_1,\ldots,\alpha_g\in F^1H^1_\dR((X,Y)/k)$, $\alpha_{g+1},\ldots,\alpha_{2g}\notin F^1\cup W_0$ and the polarization $\theta\in H^2_\dR(X/k)=\mr{Gr}^W_2H^2_\dR((X,Y)/k)$ corresponds with
$$
\theta=\alpha_1\wedge\alpha_{g+1}+\alpha_2\wedge\alpha_{g+2}+\cdots+\alpha_g\wedge\alpha_{2g},
$$
and $\alpha_0=df$ for some $f$ with $f(O)=1$ and $f(P)=0$. In particular the intersection form (in $\mr{Gr}^W_1H^1_\dR((X,Y)/\C)=H^1_\dR(X/\C)$) is given by
$$
[\langle \alpha_i,\alpha_j\rangle]=\Phi \ , \ \ \ \text{ where } \ \Phi:=\begin{pmatrix}
0 & 0 & 0 \\ 0 & 0 & I_n \\ 0 & -I_n & 0
\end{pmatrix}\in \text{Mat}_{(2n+1)\times(2n+1)}.
$$
\end{prop}

\begin{proof}
The construction of $\alpha_0$ was explained at the end of Section \ref{10may2020}. 
To construct the rest of the basis, pick any basis $\alpha_1,\ldots,\alpha_g$ of $F^1H^1_\dR((X,Y)/\C)\cong F^1H^1_\dR(X/\C)$. Since $\theta\in F^1H^2_\dR(X/\C)$ we can always write it as
$$
\theta=\alpha_1\wedge \beta_{g+1}+\alpha_2\wedge\beta_{g+2}+\cdots+\alpha_g\wedge\beta_{2g}
$$
for some $\beta_{g+i}\in H^1_\dR(X/\C)$. We claim $\alpha_1,\ldots,\alpha_g,\beta_{g+1},\ldots,\beta_{2g}$ form a basis of $H^1_\dR(X/\C)$. In fact, this follows by the equality
$$
\theta^g=g!\alpha_1\wedge\beta_{g+1}\wedge\alpha_2\wedge\beta_{g+2}\wedge\cdots\wedge\alpha_{g}\wedge\beta_{2g}.
$$
Using that $\mathsf{Tr}(\theta^g)=\deg(X)$ we get the desired basis taking $\alpha_{g+i}:=\frac{g}{\deg(X)}\beta_{g+i}$. The equality 
$$
\langle \alpha_0,\alpha_i\rangle=0,\ \ i=0, 1,2,\ldots, 2g
$$
follows from the fact $j^*\alpha_0=0\in H^1_\dR(X/\C)$, where $j:(X,\emptyset)\rightarrow (X,Y)$ is the inclusion as in \eqref{tracediagram}.
\end{proof}


\subsection{Polarized mixed Hodge structure for elliptic curves}
Using \v{C}ech cohomology with respect to the affine cover of Section \ref{computation of the relative cohomology}, we find that
\[
H^0_{\mr{dR}}(Y) \cong \{(f_0, f_1) \in \Omega^0(U_0 \cap Y) \times \Omega^0(U_1 \cap Y) \mid f_1|_{U_{01} \cap Y}= f_0|_{U_{01} \cap Y}\}.
\]
Then one may check that $(0, \frac{x-a}{x}|_{U_1 \cap Y}) \in \Omega^0(U_0 \cap Y) \times \Omega^0(U_1 \cap Y)$ represents a non-trivial element in $\mr{coker}(H^0_{\mr{dR}}(X) \to H^0_{\mr{dR}}(Y))$. We remark two things:
\begin{enumerate}
\item
The image of $[(0, \frac{x-a}{x}|_{U_1 \cap Y})]$ under the boundary map $H^0_{\mr{dR}}(Y) \to H^1_{\mr{dR}}(X, Y)$ is the same as $\omega_1 = [( (0, 0), (\mr{d}(\frac{x-a}{x}), 0), \frac{x-a}{x}|_{U_{01}})]$;
\item
The above bilinear symmetric form descends to $\mr{coker}(H^0_{\mr{dR}}(X) \to H^0_{\mr{dR}}(Y))$, which is a polarization (or the norm) on $\mr{Gr}^W_0 H^{1}_{\mr{dR}}(X, Y)$. Under this polarization, $[(0, \frac{x-a}{x}|_{U_1 \cap Y})]$ has norm $1$.
\end{enumerate}



\section{$\mathsf{T}$-space}
\label{sect5}

In this section we introduce the algebro-geometric framework where the modular vector fields defining the differential equations of quasi Jacobi forms are defined. This is the moduli space $\mathsf{T}$ of enhanced principally polarized abelian varieties with two marked points. We divide the proof of Theorem \ref{theo1} into two parts. First we show that $\mathsf{T}$ is a quasi-projective variety for any genus. In the case of elliptic curves we show that this is an affine variety and find the generators of its coordinate ring. In general, this moduli space $\mathsf{T}$ comes with a natural action of an algebraic group $\mathsf{G}$. We describe this action in the case of elliptic curves. Later, in Section \ref{sect9}, we will use this algebraic action to obtain the modularity equations satisfied by the solutions of the modular vector fields. 

\subsection{Enhanced abelian varieties with two marked points}

\begin{df}
\label{defT}
An enhanced principally polarized abelian variety with two marked points $(X,Y)$ is the data
$$
(X,Y),[\alpha_0,\alpha_1,\ldots,\alpha_{2g}]
$$
where $Y=\{O,P\},\ P\not=O$ and $\alpha_i$'s are as in  Proposition \ref{prop1}. We denote the moduli of enhanced principally polarized abelian varieties with two marked points by $\mathsf{T}$.
\end{df}
\begin{prop}
The moduli space $\mathsf{T}$ is a quasi-projective variety defined over $\Q$. 
\end{prop}

\begin{proof}
\def\Ts{{\mathsf{S}}}
The moduli  $A_g$ of principally polarized abelian varieties of dimension $g$ over a field $k$ is a quasi-projective variety  
over $\Q$, see for instance \cite{Mu1994}.
One can even construct the moduli of abelian schemes over a ring and construct the corresponding moduli 
stack, which is mainly known as Deligne-Mumford stack, however, due to the lack of motivation we avoid this and refer the reader 
to the article \cite{Fonseca} and the references therein. Let $B_g$ be the moduli space of $(A,\alpha_1,\ldots,\alpha_{2g})$ principally polarized abelian varieties equipped with a basis of $H^1_\dR(A)$ as in Proposition \ref{prop1}, this is called a symplectic-Hodge basis in \cite{Fonseca}. By \cite[Theorem 7.1]{Fonseca} this moduli space is a smooth algebraic quasi-affine variety over $\Q$ (in fact a smooth quasi-affine scheme over $\Z[\frac{1}{2}]$) with a universal family $\mathsf{U}\rightarrow B_g$. Hence $\mathsf{U}$ is the moduli of $(A,P,\alpha_1,\ldots,\alpha_{2g})$ where $P\in A$. Let $\check{\mathsf{U}}$ be the open subset of $\mathsf{U}$ corresponding to $P\neq O$ (the complement of the zero section). Then $\mathsf{T}\rightarrow\check{\mathsf{U}}$ is a line bundle over a quasi-projective variety.\end{proof}

\begin{rmk}
In the genus one case, $B_g$, $\check{\mathsf{U}}$ and so $\mathsf{T}$, are all affine varieties. In the next section we describe the coordinate ring of $\mathsf{T}$.
\end{rmk}

\subsection{Elliptic curves}
In the case of elliptic curves, after choosing the Weierstrass coordinates $x, y$, we may write $X$ in the Weierstrass format
\[
y^2 = 4x^3 - t_2x - t_3
\]
with $\Delta = 27 t_3^2 - t_2^3 \neq 0$. In these coordinates, we write $P \neq O \in Y$ as $(a, b)$. Then we get $t_3 = 4a^3 -t_2a - b^2$. 
In order to construct the basis $\alpha_0,\alpha_1,\alpha_2\in H^1_\dR(X,Y)$, first we take an arbitrary basis, for instance the basis $\omega_1,\omega_2,\omega_3$ of \eqref{24june2021}.  Note that the intersection matrix of this basis  is given by
\[
\Phi = \left(
\begin{matrix}
\langle \omega_i, \omega_j\rangle
\end{matrix}
\right) = \left(
\begin{matrix}
0 & 0 & 0 \\
0 & 0 & 1 \\
0 & -1 & 0
\end{matrix}
\right).
\]
Let us make a base change $\alpha = S \omega$ in such a way that the new basis $\alpha$ is compatible with the polarized mixed Hodge structure (as in Proposition \ref{prop1}). These restrictions imply that $S$ must be of the form
\[
\left(
\begin{matrix}
1 & 0 & 0 \\
0 &1 & 0 \\
c & t_1 &t_0
\end{matrix}
\right),
\]
where $c,t_1,t_0$ are unknown parameters. The intersection form in $\alpha_i$'s is given by
\[
[\langle \alpha_i, \alpha_j \rangle] = S [\langle \omega_i, \omega_j \rangle]S^{\mr{tr}}.
\]
Since we want to preserve the intersection matrix equal to $\Phi$, we conclude that $t_0 = 1$ and $c, t_1$ are independent parameters. Therefore
$$
{\sf{T}}:=\mr{Spec}\ \mb{C}[a, b, c, t_1, t_2, \frac{1}{\Delta}]
$$ 
and over $\sf{T}$,
we have the universal family
\begin{equation} \nonumber
\begin{split}
&X: y^2 = 4x^3 - t_2x -t_3,(a, b), \\
&\mr{d}\left(\frac{x-a}{x}\right) \ , \ \  \frac{\mr{d}x}{y} \ , \ \   \left(c-\frac{b}{2a}\right)\mr{d}\left(\frac{x-a}{x}\right)+t_1\frac{\mr{d}x}{y} + \frac{x\mr{d}x}{y}-\mr{d}\left(\frac{y}{2x}\right).
\end{split}
\end{equation} 

\begin{rmk}
For a general construction of this enhanced moduli space $\mathsf{T}$ for other families of varieties see \cite[Theorem 3.5]{M2020}.
\end{rmk}

\subsection{Algebraic group}
\label{alggrp}
We define the algebraic group $\mathsf{G}$ to be the automorphism group ${\rm Aut}(H^*_{dR}(X_0, Y_0), W_0^*,F_0^*, \cup, \theta_0)$ of a fixed enhanced abelian variety $X_0$, where $\theta_0$ is the polarization of $X_0$, see \cite[Section 3.3, Section 11.4]{M2020}. The notation ${\rm Aut}(H^*_{dR}(X_0, Y_0), W_0^*,F_0^*, \cup, \theta_0)$ means we are considering the automorphisms of $H^*_{dR}(X_0, Y_0)$ which respects the polarized mixed Hodge structures and preserves the intersection matrices. In the case of abelian varieties with two marked points, it can be computed explicitly
$$
\mathsf{G}=\left\{\begin{pmatrix}1 & 0 & v \\ 0 & k & k' \\ 0 & 0 & k^{-\mathsf{tr}} \end{pmatrix}\in \text{GL}_{2g+1}(\C) : \ kk'^{\mathsf{tr}}=k'k^\mathsf{tr} \ , \ v\in \C^g\right\},
$$
where the form of $\mathsf{g}\in\mathsf{G}$ is derived from the fact that it respects the weight and Hodge filtration, and hence $\mathsf{g}^{12}=\mathsf{g}^{21}=\mathsf{g}^{31}=\mathsf{g}^{32}=0$, also $\mathsf{g}^{\mathsf{tr}}\Phi\mathsf{g}=\Phi$ and the entry equal to $1$ is given in order to preserve the polarization on $Gr^W_0H^1_\dR(X,Y)$. The algebraic group $\mathsf{G}$ is of dimension $\frac{3g(g+1)}{2}$ and it acts on $\sf{T}$ by change of basis of the de Rham cohomology group $H^1_{dR}(X, Y)$.

\subsection{The algebraic group for elliptic curves}
In the case of elliptic curves  
$\mathsf{G}$ acts on the de Rham classes $\alpha_i, i = 1,2,3$ in the following way
\[
(\alpha_0, \alpha_1,\alpha_2) \bullet \mathsf{g} = (\alpha_0, k \alpha_1, v \alpha_0+k' \alpha_1+k^{-1}\alpha_2)
\]
for $\mathsf{g} = \left(
\begin{matrix}
1 & 0 & v \\
0 & k & k' \\
0 & 0 & k^{-1}
\end{matrix}
\right) \in \mathsf{G}$.
This action induces an action on $\mathsf{T}$, which we can describe explicitly in terms of its parameters. Consider the parameter $(a,b,c,t_1,t_2)\in\mathsf{T}$ corresponding to the enhanced elliptic curve
\[
X: y^2 = 4x^3 -t_2x -t_3 = 4x^3 - t_2x + b^2 - 4a^3 +t_2a, 
\]
together with the marked point $P=(a,b)$, and the frame 
\[
(\alpha_0,\alpha_1,\alpha_2)=(\omega_1,\omega_2,c\cdot\omega_1+t_1\cdot\omega_2+\omega_3),
\]
where $\omega_1,\omega_2,\omega_3\in H^1_\dR(X,Y)$ is the basis given in Proposition \ref{basis of relative cohomology}, this is
\[
\omega_1=\mr{d}\left(\frac{x-a}{x}\right) \ , \ \  \omega_2=\frac{\mr{d}x}{y} \ , \ \   \omega_3=\frac{x\mr{d}x}{y}-\mr{d}\left(\frac{y}{2x}\right)-\frac{b}{2a}\mr{d}\left(\frac{x-a}{x}\right).
\]
Applying $\mathsf{g}$ 
to the point $(a,b,c,t_1,t_2)$ above we get the same elliptic curve, with the same marked points, but now the frame is
$$
(\alpha_0',\alpha_1',\alpha_2')=(\omega_1,k\cdot\omega_2,(v+k^{-1}c)\cdot\omega_1+(k'+k^{-1}t_1)\cdot\omega_2+k^{-1}\cdot\omega_3).
$$
In order to see which is the point in $\mathsf{T}$ corresponding to it, we consider the following change of coordinates
\[
\varphi: \mb{A}^2 \to \mb{A}^2, (x, y) \to (k^2x, k^3 y).
\]
Under this map, the resulting enhanced elliptic curve $(\varphi^{-1}(X), \varphi^{-1}(Y))$, $\varphi^*[\alpha_0',\alpha_1',\alpha_2']$ is given by
\[
y^2=4x^3-k^{-4}t_2x-k^{-6}t_3,
\]
with the marked point $\varphi^{-1}(P)=(k^{-2}a,k^{-3}b)$, and the frame
\[
\varphi^*\alpha_0'=\mr{d}\left(\frac{x-k^{-2}a}{x}\right) \ , \ \  \varphi^*\alpha_1'=\frac{\mr{d}x}{y} \ , 
\]
\vspace{-3mm}
\[
\varphi^*\alpha_2'=(v+k^{-1}c)\varphi^*\alpha_0'+(k^{-1}k'+k^{-2}t_1)\varphi^*\alpha_1'+\frac{x\mr{d}x}{y}-\mr{d}\left(\frac{y}{2x}\right)-\frac{k^{-3}b}{2k^{-2}a}\mr{d}\left(\frac{x-k^{-2}a}{x}\right).
\]
Then we derive the action of $\mathsf{g} \in \mathsf{G}$ on $\mathsf{T}$ via $S \bullet \mathsf{g} = \mathsf{g}^{\mr{tr}} \cdot S \cdot \mr{diag}(1, k^{-1}, k)$, this is 
\[
t = (a, b, c, t_1, t_2) \to t \bullet \mathsf{g} = (k^{-2}a, k^{-3}b, v+k^{-1}c, k^{-1}k' + k^{-2}t_1, k^{-4}t_2).
\]

\section{Gauss-Manin connection}
\label{section6}
This section is devoted to the computation of the Gauss-Manin connection. In order to do this, we use the fundamental property of the Gauss-Manin connection in relation to integrals. After using this analytical interpretation we provide another interpretation of the same computation from the algebraic description of the Gauss-Manin connection. 
\subsection{Gauss-Manin connection and integrals}
Recall the following exact sequence of a pair $(X,Y)$
\[
0 \to H^0_{\mr{dR}}(X) \to H^{0}_{\mr{dR}}(Y) \to H^{1}_{\mr{dR}}(X, Y) \to H^1_{\mr{dR}}(X) \to 0.
\]
When the pair $(X, Y)$ varies, the above exact sequence can be viewed as an exact sequence of local systems.
\begin{df}
The Gauss-Manin connection on $H^{1}_{\mr{dR}}(X, Y)$ is a connection (i.e. $\mb{C}$-linear with the Leibniz rule) on $H^{1}_{\mr{dR}}(X, Y)$
\[
\nabla: H^{1}_{\mr{dR}}(X, Y) \to \Omega^1_T \otimes H^{1}_{\mr{dR}}(X, Y)
\]
satisfies that $\nabla$ is compatible with the above exact sequence, i.e.,
\begin{equation}
\begin{split}
0 &\to (H^0_{\mr{dR}}(X), \nabla) \to (H^{0}_{\mr{dR}}(Y), \nabla) \\
&\to (H^{1}_{\mr{dR}}(X, Y), \nabla) \to (H^1_{\mr{dR}}(X), \nabla^{}) \to 0.
\end{split}
\end{equation}
Similar to the non-relative case, the relative Gauss-Manin connection is uniquely determined by its fundamental property 
\begin{equation}
 d\left(  \int_{\delta}\omega\right)=\int_{\delta}\nabla\omega,\ \ \ \ \ \delta\in H_1(X,Y),\ \ \omega\in H^1_\dR(X,Y)
\end{equation}
where the second integration occurs only in the $H^1_\dR(X,Y)$ piece, in other words if
$$
\nabla\omega=\sum_{i\in I}df_i\otimes \eta_i,
$$
then
$$
\int_\delta\nabla\omega:=\sum_{i\in I}df_i\cdot\int_\delta \eta_i.
$$
\end{df}
\subsection{Elliptic curves}
Recall that $\alpha_1 = \mr{d}(\frac{x-a}{x}), \alpha_2 = \frac{\mr{d}x}{y}$ and $\alpha_3 = (c-\frac{b}{2a})\mr{d}(\frac{x-a}{x})+t_1\frac{\mr{d}x}{y} + \frac{x\mr{d}x}{y}-\mr{d}(\frac{y}{2x})$. Recall $t_3 = 4a^3-b^2-at_2$ and hence $\mr{d}t = (12a^2-t_2)\mr{d}a - 2b \mr{d}b -a\mr{d}t_2$. We also set $\alpha = 3t_3 \mr{d}t_2-2t_2\mr{d}t_3$.
\begin{prop}
The Gauss-Manin connection of the family of elliptic curves
\[
(X: y^2 = 4x^3 - t_2x - t_3, (a, b), \alpha_1, \alpha_2, \alpha_3)
\]
is given as follows:
\begin{equation}
\nabla \left(\begin{matrix}
\alpha_1 \\
\alpha_2 \\
\alpha_3
\end{matrix} \right) = \left(\begin{matrix}
0 & 0 & 0 \\
A_{21} & A_{22} & A_{23} \\
A_{31} & A_{32} & A_{33} 
\end{matrix} \right) \otimes \left(\begin{matrix}
\alpha_1 \\
\alpha_2 \\
\alpha_3
\end{matrix} \right)
\end{equation}
where
\begin{small}
\begin{eqnarray} \nonumber
A_{21} &=&  g_1\mr{d}t_2+g_2 \mr{d}t_3 - \left(c-\frac{b}{2a}\right)\frac{3\alpha}{2}-\frac{\Delta\mr{d}a}{b}, \\ \nonumber
A_{22} &=&  -\frac{3t_1\alpha}{2}-\frac{\mr{d}\Delta}{12},  \\ \nonumber
A_{23} &=& \frac{3\alpha}{2} ,  \\ \nonumber
A_{31} &=&  (t_1g_1+g_3)\mr{d}t_2+(t_1g_2+g_1)\mr{d}t_3 -\frac{(a+t_1)\Delta \mr{d}a}{b}- \left(c-\frac{b}{2a}\right)\left(\frac{3t_1\alpha}{2}+\frac{\mr{d}\Delta}{12}\right)+\Delta \mr{d}c, \\ \nonumber
A_{32} &=&  \Delta \mr{d}t_1-\frac{t_1\mr{d}\Delta}{6}-\left(\frac{3t_1^2}{2}+\frac{t_2}{8}\right)\alpha,  \\ \nonumber
A_{33} &=&  \frac{3t_1\alpha}{2}+\frac{\mr{d}\Delta}{12}, 
\end{eqnarray}
\end{small}
and 
\begin{eqnarray} \nonumber
g_1 &=& \frac{-2a^2t_2^2+3at_2t_3+9t_3^2}{4ab}, \\ \nonumber
g_2 &=& \frac{18a^2t_3-at_2^2-3t_2t_3}{2ab},  \\ \nonumber
g_3 &=& \frac{6a^2t_2t_3+(18t_3^2-t_2^3)a-t_2^2t_3}{8ab}.
\end{eqnarray}
\end{prop}
\begin{proof}
In the absolute case, the proof is a classical calculation, see \cite[Proposition 3.1]{M2012} for example. For the relative case, we mention a theorem of Fuchs (see \cite[Theorem 1.1]{Manin}) from the historical point view, whose proof contains the computation of relative Gauss-Manin connection. Let us compute the relative Gauss-Manin connection in the basis
$$
\omega_1=\mr{d}\left(\frac{x-a}{x}\right),\  \omega_2=\frac{\mr{d}x}{y},\ \omega_3=\frac{x\mr{d}x}{y}-\mr{d}\left(\frac{y}{2x}\right).
$$
For this we will use the second's author computation of the absolute Gauss-Manin connection in the basis $\omega_2, \omega_3$ (see \cite[Proposition 3.1]{M2012})
\begin{equation}
\label{absGMC}
\def\arraystretch{2}
\begin{pmatrix}
\mr{d}\left({\int}\omega_2\right)\\
\mr{d}\left({\int} \omega_3 \right)\
\end{pmatrix}
=
\begin{pmatrix}
{-\frac{1}{12}\frac{\mr{d}\Delta}{\Delta}} &
{\frac{3}{2}\frac{\alpha}{\Delta}} \\
{ -\frac{1}{8}t_2\frac{\alpha}{\Delta}} &
{\frac{1}{12}\frac{\mr{d}\Delta}{\Delta}}
\end{pmatrix}
\begin{pmatrix}
{\int}\omega_2\\
{\int} \omega_3  
\end{pmatrix}
\end{equation}
where $\mr{d}$ is the differential with respect to $t_2, t_3$ and the integration is taken along any closed path inside $X$. The main relation between the absolute Gauss-Manin connection and the relative one can be obtained as follows: Assume that we know the following computation for the family $y^2=p(x)$
\begin{equation}
\label{affrelGMC}
\def\arraystretch{2}
\begin{pmatrix}
\mr{d}\left({\int}\frac{dx}{y}\right)\\
\mr{d}\left({\int} \frac{x\mr{d}x}{y} \right)\
\end{pmatrix}
=
\begin{pmatrix}
{-\frac{1}{12}\frac{\mr{d}\Delta}{\Delta}} &
{\frac{3}{2}\frac{\alpha}{\Delta}} \\
{ -\frac{1}{8}t_2\frac{\alpha}{\Delta}} &
{\frac{1}{12}\frac{\mr{d}\Delta}{\Delta}}
\end{pmatrix}
\begin{pmatrix}
{\int}\frac{\mr{d}x}{y}\\
{\int} \frac{x\mr{d}x}{y}  
\end{pmatrix}
+\begin{pmatrix}
\left(\int \mr{d}h_1\right)\mr{d}t_2+\left(\int \mr{d}h_2\right)\mr{d}t_3\\
\left(\int \mr{d}h_3\right)\mr{d}t_2+\left(\int \mr{d}h_4\right)\mr{d}t_3
\end{pmatrix}
\end{equation}
where $\mr{d}$ is the differential with respect to $t_2,t_3$ and the integration is over any path in $X$ minus $O$.
We write this equality with the correction of $\frac{x\mr{d}x}{y}$, taking differential with respect to parameters $t_2,t_3,a$.
After this we can take  the path of integration from $O$ to $P$ and obtain
\begin{equation}
 \label{21apr2021}
\def\arraystretch{2}
\begin{pmatrix}
\mr{d}\left({\int}\frac{\mr{d}x}{y}\right)\\
\mr{d}\left({\int}\left(\frac{x\mr{d}x}{y}-\mr{d}(\frac{y}{2x})\right)\right)\
\end{pmatrix}
=
\begin{pmatrix}
{-\frac{1}{12}\frac{\mr{d}\Delta}{\Delta}} &
{\frac{3}{2}\frac{\alpha}{\Delta}} \\
{ -\frac{1}{8}t_2\frac{\alpha}{\Delta}} &
{\frac{1}{12}\frac{\mr{d}\Delta}{\Delta}}
\end{pmatrix}
\begin{pmatrix}
{\int}\frac{\mr{d}x}{y}\\
{\int}  \left(\frac{x\mr{d}x}{y}-\mr{d}(\frac{y}{2x})\right)  \
\end{pmatrix}
+\begin{pmatrix}
X \\ Y
\end{pmatrix}
\end{equation}
\begin{equation}
\begin{smallmatrix}
X = & \left(\int \mr{d}h_1\right)\mr{d}t_2+\left(\int \mr{d}h_2\right)\mr{d}t_3  + \left({\int} \mr{d}(\frac{y}{2x})\right )
\frac{3\alpha}{2\Delta}
+\frac{\mr{d}a}{b},\\
Y = & \left(\int \mr{d}(h_3+\frac{1}{4y})\right)\mr{d}t_2+\left(\int \mr{d}(h_4+\frac{1}{4xy})\right)\mr{d}t_3 + \left({\int} \mr{d}(\frac{y}{2x})\right)\frac{\mr{d}\Delta}{12\Delta} +    \frac{a\mr{d}a}{b}-\mr{d}(\frac{b}{2a})-\frac{\mr{d}t_3+a\mr{d}t_2}{4ab}.
\end{smallmatrix}
\end{equation}
This relation tells us the relative Gauss-Manin connection in the basis $\omega_1$, $\omega_2$, $\omega_3$ is:
\begin{equation}
\label{relGMC}
\def\arraystretch{2}
\begin{pmatrix}
\nabla(\omega_1)\\
\nabla(\omega_2)\\
\nabla(\omega_3)
\end{pmatrix}
=
\begin{pmatrix}
0 &  0 & 0 \\
-X  & {-\frac{1}{12}\frac{d\Delta}{\Delta}} &
{\frac{3}{2}\frac{\alpha}{\Delta}} \\
-Y  & { -\frac{1}{8}t_2\frac{\alpha}{\Delta}} &
{\frac{1}{12}\frac{d\Delta}{\Delta}}
\end{pmatrix}
\begin{pmatrix}
\omega_1\\
\omega_2\\
\omega_3
\end{pmatrix}.
\end{equation}
In order to compute $h_1$, $h_2$, $h_3$ and $h_4$ appearing in \eqref{affrelGMC}, we compute the Gauss-Manin connection of $\omega_2$ and $\omega_3$. To do this we have to write each
$$
\mr{d}\omega_i=\sum_{a}\alpha_a\wedge\beta_a\in\Omega_T^1\wedge\Omega_{E/T}^1, 
$$
where $\mr{d}$ is the differential with respect to $x$, $y$, $t_2$, $t_3$, and $a$.
Then 
$$
\nabla\omega_i=\sum_a\alpha_a\otimes\beta_a.
$$
Note first that since $\omega_1=\mr{d}\left(\frac{x-a}{x}\right)$, then $\nabla\omega_1=0$ and the first row of \eqref{relGMC} follows. Using the equation 
$$
2y\mr{d}y=\mr{d}p
$$
we get
$$
\mr{d}\omega_2=\frac{(x\mr{d}t_2+\mr{d}t_3)\wedge \mr{d}x}{2y^3} \ , \ \ \ \ \mr{d}\omega_3=\frac{(x^2\mr{d}t_2+x\mr{d}t_3)\wedge \mr{d}x}{y}
$$
and so
$$
\nabla\omega_2=\mr{d}t_2\otimes \frac{x\mr{d}x}{2y^3}+\mr{d}t_3\otimes \frac{\mr{d}x}{2y^3} \ , \ \ \ \nabla\omega_3=\mr{d}t_2\otimes \frac{x^2\mr{d}x}{2y^3}+\mr{d}t_3\otimes \frac{x\mr{d}x}{2y^3}.
$$
Following the method described in \cite[Proposition 2.2]{M2012}, we can obtain the following identities:
$$
\frac{x^2\mr{d}x}{y}=\mr{d}\left(\frac{y}{6}\right)+\frac{t_2}{12}\frac{\mr{d}x}{y},
$$
$$
\frac{x^3\mr{d}x}{y}=\mr{d}\left(\frac{xy}{10}\right)+\frac{3t_2}{20}\frac{x\mr{d}x}{y}+\frac{t_3}{10}\frac{\mr{d}x}{y},
$$
$$
\frac{x^4\mr{d}x}{y}=\mr{d}\left(\frac{x^2y}{14}+\frac{5t_2y}{168}\right)+\frac{t_3}{7}\frac{x\mr{d}x}{y}+\frac{5t_2^2}{336}\frac{\mr{d}x}{y},
$$
$$
\frac{x^5\mr{d}x}{y}=\mr{d}\left(\frac{x^3y}{18}+\frac{7t_2xy}{360}+\frac{t_3y}{36}\right)+\frac{7t_2^2}{240}\frac{x\mr{d}x}{y}+\frac{t_2t_3}{30}\frac{\mr{d}x}{y}.
$$
Using the following identity
$$
\Delta=-(12x^2-t_2)A+pB
$$
where 
$$
A=-36x^4+15t_2x^2-t_2^2,
$$
$$
B=-108x^3+27t_2x-27t_3,
$$
together with the previous ones, we obtain
$$
\frac{\mr{d}x}{2y^3}=\frac{\mr{d}x}{2py}=\frac{1}{\Delta}\frac{(-(12x^2-t_2)A+pB)\mr{d}x}{2py}=\frac{1}{\Delta}\left(-\frac{A\mr{d}y}{y^2}+\frac{B}{2}\frac{\mr{d}x}{y}\right)
$$
$$
=\frac{1}{\Delta}\left(\mr{d}\left(\frac{A}{y}\right)+\left(\frac{B}{2}-A'\right)\frac{\mr{d}x}{y}\right)
$$
$$
=\mr{d}\left(\frac{A}{\Delta y}\right)+\frac{1}{\Delta}\left(90x^3-\frac{33t_2}{2}x-\frac{27t_3}{2}\right)\frac{\mr{d}x}{y}
$$
$$
=\mr{d}\left(\frac{A}{\Delta y}\right)+\frac{1}{\Delta}\left(\mr{d}(9xy)-3t_2\frac{x\mr{d}x}{y}-\frac{9t_3}{2}\frac{\mr{d}x}{y}\right)
$$
$$
=\mr{d}\left(\frac{A}{\Delta y}+\frac{9xy}{\Delta}\right)-\frac{1}{\Delta}\left(3t_2\frac{x\mr{d}x}{y}+\frac{9t_3}{2}\frac{\mr{d}x}{y}\right).
$$
Hence 
$$
h_2=\frac{A}{\Delta y}+\frac{9xy}{\Delta}.
$$
Similarly we compute
$$
\frac{x\mr{d}x}{2y^3}=\mr{d}\left(\frac{Ax}{\Delta y}+\frac{9x^2y}{\Delta}-\frac{3t_2y}{2\Delta}\right)+\frac{1}{\Delta}\left(\frac{9t_3}{2}\frac{x\mr{d}x}{y}+\frac{t_2^2}{4}\frac{\mr{d}x}{y}\right),
$$
$$
\frac{x^2\mr{d}x}{2y^3}
=\mr{d}\left(\frac{Ax^2}{\Delta y}+\frac{9x^3y}{\Delta}-\frac{3t_2xy}{2\Delta}+\frac{9t_3y}{4\Delta}\right)-\frac{t_2^2}{4\Delta}\frac{x\mr{d}x}{y}-\frac{3t_2t_3}{8\Delta}\frac{\mr{d}x}{y}
$$
and so we obtain
$$
h_1=h_4=\frac{Ax}{\Delta y}+\frac{9x^2y}{\Delta}-\frac{3t_2y}{2\Delta} \ , \ \ \ h_3=\frac{Ax^2}{\Delta y}+\frac{9x^3y}{\Delta}-\frac{3t_2xy}{2\Delta}+\frac{9t_3y}{4\Delta}.
$$
Replacing these functions $h_i, i=1,2,3,4$ in \eqref{relGMC}, we obtain that the Gauss-Manin matrix with respect to $\omega_1$, $\omega_2$ and $\omega_3$ is
\[
B =  \frac{1}{\Delta}\left(\begin{matrix}
0 & 0 & 0\\
B_{21} & -\frac{\mr{d}\Delta}{12} & \frac{3\alpha}{2} \\
B_{31} & -\frac{t_2 \alpha}{8} & \frac{\mr{d} \Delta}{12}
\end{matrix}\right),
\]
where $\alpha = 3t_3 \mr{d}t_2 - 2t_2 \mr{d}t_3, $
\begin{eqnarray} \label{B-matrix}
B_{21} &=& g_1 \mr{d}t_2 + g_2 \mr{d}t_3 - \frac{\Delta\mr{d}a}{b}, \\ \nonumber
B_{31} &=& g_3 \mr{d}t_2 + g_1 \mr{d}t_3 - \frac{a\Delta\mr{d}a}{b}+\Delta \mr{d}\left(\frac{b}{2a}\right)
\end{eqnarray}
and $g_i$ is defined as above.
For simplicity, we write $B$ as
\[B = \frac{1}{\Delta}\left(\begin{matrix}
0 & 0 & 0\\
B_{21} & B_{22} & B_{23} \\
B_{31} & B_{32} & B_{33}
\end{matrix}\right).\]
We let \[A = \frac{1}{\Delta}\left(\begin{matrix}
0 & 0 & 0\\
A_{21} & A_{22} & A_{23} \\
A_{31} & A_{32} & A_{33}
\end{matrix}\right),\]
be the Gauss-Manin matrix under the basis $\omega_1, \omega_2, (c-\frac{b}{2a})\omega_1+t_1\omega_2 + \omega_3$. Then
\begin{small}
\begin{equation} \nonumber
\begin{split}
&\nabla\left(\frac{\mr{d}x}{y}\right) = \frac{1}{\Delta}\left(B_{21} \mr{d}\left(\frac{x-a}{x}\right) + B_{22} \frac{\mr{d}x}{y} + B_{23}\left(\frac{x\mr{d}x}{y}-\mr{d}\left(\frac{y}{2x}\right)\right)\right) \\
&= \frac{1}{\Delta}\left(A_{21} \mr{d}\left(\frac{x-a}{x}\right) + A_{22} \frac{\mr{d}x}{y} + A_{23}\left(\left(c-\frac{b}{2a}\right)\mr{d}\left(\frac{x-a}{x}\right)+t_1\frac{\mr{d}x}{y} + \frac{x\mr{d}x}{y}-\mr{d}\left(\frac{y}{2x}\right)\right)\right).
\end{split}
\end{equation}
\end{small}
Therefore we get
\begin{equation} \label{A-matrix1}
A_{21} = B_{21} -\left(c-\frac{b}{2a}\right) B_{23} \ , \ \ \  A_{22} = B_{22} - dB_{23} \ , \ \ \ A_{23} = B_{23}.
\end{equation}
Similarly we have
\begin{eqnarray} \nonumber
&A_{31}& = t_1B_{21}+B_{31} -  \left(c-\frac{b}{2a}\right)(t_1B_{23}+B_{33})+\Delta \mr{d}\left(c-\frac{b}{2a}\right),  \\  \nonumber
&A_{32}& = t_1B_{22}+B_{32} - t_1(t_1B_{23}+B_{33}) + \Delta \mr{d}t_1 \ , \ \ \ \  A_{33} = t_1B_{23} + B_{33}.
\end{eqnarray}
This can be also proved by using the formula 
\[A = \mr{d}S \cdot S^{-1} + S \cdot B \cdot S^{-1},\]
where $S = \left( \begin{matrix}
1 & 0 & 0 \\
0 & 1 & 0 \\
c-\frac{b}{2a} & t_1 & 1 
\end{matrix}\right).$
Finally we get the desired expressions.
\end{proof}

\subsection{Algebraic Gauss-Manin connection for elliptic curves}
\label{subsectalgGMell}
In the case of elliptic curves, we have another way to describe the relative Gauss-Manin connection by means of the isomorphism \eqref{22june2021} and the canonical Gauss-Manin connection on $R^1\pi_*\Omega_{(X,Y)/\mathsf{T}}^\bullet$.

\begin{prop}
\label{isoeq2.6familia}
Let $\pi:\mathsf{X}\rightarrow \mathsf{T}$ be the family of elliptic curves with two marked points. Let $\mathsf{Y}\subseteq\mathsf{X}$ be the subvariety of $\mathsf{X}$ such that $\mathsf{Y}_t=\{O,P\}$ for every $t\in \mathsf{T}$. Let $V\subseteq \mathsf{T}$ be an affine open set,  $\mathsf{U}\subseteq\mathsf{X}_V$ be an affine open set such that $\mathsf{Y}_V\subseteq \mathsf{U}$, $\mathsf{U}\rightarrow V$ is an affine locally trivial fibration, and there exists some $f\in \Gamma(\mathcal{O}_{\mathsf{U}_V})$ with $f(O)=1$ and $f(P)=0$ for all $t\in V$. Then the isomorphism of Proposition \ref{isoeq2.6} extends to the family
$$
H^1_\dR((\mathsf{X},\mathsf{Y})/V)=\uhp^1(\mathsf{X}_V, \Omega_{(\mathsf{X},\mathsf{Y})/V}^\bullet)\cong H^1(\Gamma(\hat{\Omega}_{\mathsf{U}/V}^\bullet)).
$$
\end{prop}

\begin{proof}
Just note that $R^1\pi_*\Omega_{(\mathsf{X},\mathsf{Y})/V}^\bullet$ is a vector bundle of rank 3 over $V$ and the same holds for $R^1\pi_*\hat{\Omega}_{\mathsf{U}/V}^\bullet$. Under the hypothesis of the theorem we have a well defined morphism of bundles
$$
R^1\pi_*\Omega_{(\mathsf{X},\mathsf{Y})/V}^\bullet\rightarrow R^1\pi_*\hat{\Omega}_{\mathsf{U}/V}^\bullet
$$
which is in fact an isomorphism, since it is an isomorphism on each fiber. Therefore it induces an isomorphism in global sections, and since $\hat{\Omega}_{\mathsf{U}/V}^\bullet$ is acyclic, we get the result.
\end{proof}

\begin{rmk}
If we denote $\pi_0:\mathsf{X}_0\rightarrow \mathsf{T}_0$ the usual family of elliptic curves, we have a similar result for the absolute de Rham cohomology. In fact, given $V_0\subseteq\mathsf{T}_0$ an affine open set and $\mathsf{U}_0\subseteq\mathsf{X}_0$ an affine open set such that $\mathsf{U}_0\rightarrow V_0$ is an affine locally trivial fibration, then we have an isomorphism
$$
H^1_\dR(\mathsf{X}_0/V_0)=\uhp^1((\mathsf{X}_0)_{V_0},\Omega_{\mathsf{X}_0/V_0}^\bullet)\cong H^1(\Gamma(\bar{\Omega}_{\mathsf{U}_0/V_0}^\bullet)),
$$
where now $\bar{\Omega}_{\mathsf{U}_0/V_0}^\bullet$ is the relative (to the family) version of the complex
$$
\bar{\Omega}_U^\bullet: \hspace{0.5cm} i_*\mathcal{O}_U\xrightarrow{d}i_*\Omega_U^1\xrightarrow{Res} \C_{X- U}\rightarrow 0 \ .
$$
This isomorphism is explicitly described as follows: Let $X$ be an elliptic curve, $U\subseteq X$ be an affine open set and let $V=X- \{Q\}$, such that $U\cup V=X$. Then a proof similar to Proposition \ref{isoeq2.6} shows that the map
$$
\omega=(\omega_0,\omega_1,f_{01})\in \uhp^1(\mathcal{U},\Omega_X^\bullet)
$$
$$
\mapsto \omega_0\in H^1(\Gamma(\bar{\Omega}_U^\bullet))=\frac{\omega\in\Gamma(\Omega_U^1)\text{ without residues on }X- U}{\text{exact forms $df$ with }f\in\Gamma(\mathcal{O}_U)}
$$
is an isomorphism. Moreover, in \cite[Section 3.3]{M2012} is described how one can compute the Gauss-Manin connection directly on $H^1(\Gamma(\bar{\Omega}_{\mathsf{U}_0/V_0}^\bullet))$ just considering for any $\alpha\in H^1(\Gamma(\bar{\Omega}_{\mathsf{U}_0/V_0}^\bullet))$
$$
d\alpha=\sum_{i}   dt_i\wedge\beta_i
$$
and then $\nabla(\alpha)=\sum_i dt_i\otimes \beta_i\in \Omega^1_{\mathsf{T}_0}(V_0)\otimes H^1(\Gamma(\bar{\Omega}_{\mathsf{U}_0/V_0}^\bullet))$. On the relative case we have the following result, which turns out to be equivalent to our computation of the previous section.
\end{rmk}

\begin{prop}
With the same notation of Proposition \ref{isoeq2.6familia}. Let $\alpha\in H^1(\Gamma(\hat{\Omega}_{\mathsf{U}/V}))$. If
$$
d\alpha=\sum_i  dt_i\wedge\beta_i
$$
then
$$
\nabla(\alpha)=(\alpha|_O-\alpha|_P)\otimes df+\sum_i dt_i\otimes \beta_i\in\Omega_{\mathsf{T}}^1(V)\otimes H^1(\Gamma(\hat{\Omega}_{\mathsf{U}/V}^\bullet)).
$$
\end{prop}

\begin{proof}
Let $\omega=((\omega_0,f_0),(\omega_1,f_1),f_{01})\in \uhp^1(\mathsf{X}_V,\Omega_{(\mathsf{X},\mathsf{Y})/V}^\bullet)$ such that
$$
\alpha=\omega_1+(f_1(P)-f_1(O))df.
$$
In order to compute $\nabla(\omega)$ we compute
$$
D\omega=((d\omega_0,\omega_0|_{Y\cap U_0}-df_0), (d\omega_1,\omega_1|_{Y\cap U_1}-df_1), (\omega_1-\omega_0|_{U_{01}}-df_{01}, f_1-f_0-f_{01}|_{Y\cap U_{01}}))
$$
$$
=\sum_i ((\beta^i_0,h^i_0),(\beta^i_1,h^i_1), h^i_{01})\wedge dt_i,
$$
and so 
$$
\nabla(\omega)=-\sum_i dt_i\otimes ((\beta^i_0,h^i_0),(\beta^i_1,h^i_1),h^i_{01}).
$$
Using the isomorphism we get
$$
\nabla(\alpha)=-\sum_i dt_i\otimes (\beta^i_1+(h^i_1(P)-h^i_1(O))df)=\sum_i dt_i\otimes \beta_i
$$
with $\beta_i=-\beta^i_1-(h^i_1(P)-h^i_1(O))df$. The result follows once we note that
$$
d\alpha=d\omega_1+(df_1(P)-df_1(O))\wedge df
$$
$$
=\sum_i  \beta^i_1\wedge dt_i-(h^i_1(P)-h^i_1(O))dt_i\wedge df+(\omega_1|_{P\cap U_0}-\omega_1|_{O\cap U_0})\wedge df
$$
$$
=(\alpha|_P-\alpha|_O)\wedge df+\sum_i dt_i\wedge\beta_i.
$$
Where in the last equality we used that $df|_Y=0$.
\end{proof}

\section{Modular vector fields}
\label{sect7}

Using the Gauss-Manin computations of the previous section we will describe the modular vector fields over the moduli space $\mathsf{T}$ of enhanced elliptic curves with two marked points. Relying in our algebraic description of $\mathsf{T}$, we solve explicitly this algebraic problem. On the other hand, for the general case of enhanced abelian varieties we will show the existence of the modular vector fields in Section \ref{gpd2021} using the period map, hence by transcendental methods.

\subsection{Proof of Theorem \ref{23.06.2021}}
The proof is based on explicit calculations. Let $g_i$ as the previous section. Then we may decompose the GM matrix under the basis $\alpha_i$ as
\begin{equation}
 A = \frac{1}{\Delta}(A_a \otimes \mr{d}a+ A_b \otimes \mr{d}b + A_c \otimes \mr{d}c +  A_{t_1} \otimes \mr{d}t_1+ A_{t_2} \otimes \mr{d}t_2),
\end{equation}
where 
\begin{tiny}
$$
A_a = 
\left(\begin{array}{cc}
0 & 0  \\
(g_2+3t_2(c-\frac{b}{2a}))(12a^2-t_2)-\frac{\Delta}{b} &(3t_1t_2-\frac{9t_3}{2})(12a^2-t_2)\\
(t_1g_2+g_1+(3t_1t_2-\frac{9t_3}{2})(c-\frac{b}{2a}))(12a^2-t_2)-\frac{(a+t_1)\Delta}{b}&(-9t_1t_3+3t_1^2t_2+\frac{t_2^2}{4})(12a^2-t_2) 
\end{array} \right.
$$
$$
\left. \begin{array}{c}
0\\
3t_2^2-36a^2t_2 \\
(-3t_1t_2+\frac{9t_3}{2})(12a^2-t_2)
\end{array}\right), 
$$
$$
A_b = 
\left(\begin{array}{ccc}
0 & 0 &0  \\
(g_2+3t_2(c-\frac{b}{2a}))(-2b) & (3t_1t_2-\frac{9t_3}{2})(-2b) & 6bt_2\\
(t_1g_2+g_1+(3t_1t_2-\frac{9t_3}{2})(c-\frac{b}{2a}))(-2b) & (-9t_1t_3+3t_1^2t_2+\frac{t_2^2}{4})(-2b) & (-3t_1t_2+\frac{9t_3}{2})(-2b)
\end{array} \right), 
$$
$$
A_{c} = \left(\begin{array}{ccc}
0 & 0 & 0\\
0 & 0 & 0\\
\Delta & 0 & 0
\end{array} \right),
A_{t_1} = \left(\begin{array}{ccc}
0 & 0 & 0\\
0 & 0 & 0\\
0 & \Delta & 0
\end{array} \right).
$$
$$
A_{t_2} = 
\left(\begin{array}{c}
0   \\
g_1-\frac{9t_3}{2}(c-\frac{b}{2a})+(g_2+3t_2(c-\frac{b}{2a}))(-a) \\
t_1g_1+g_3 + (\frac{t_2^2}{4}- \frac{9t_1t_3}{2})(c-\frac{b}{2a})+(t_1g_2+g_1+(3t_1t_2-\frac{9t_3}{2})(c-\frac{b}{2a}))(-a)
\end{array} \right.
$$
$$
\left. \begin{array}{cc}
0 & 0\\
-\frac{9t_1t_3}{2}+\frac{t_2^2}{4}+(3t_1t_2-\frac{9t_3}{2})(-a) & \frac{9t_3}{2}+6at_2 \\
\frac{t_1t_2^2}{2}-\frac{9t_1^2t_3}{2}-\frac{3t_2t_3}{8}+9at_1t_3-3at_1^2t_2-\frac{at_2^2}{4} & \frac{9t_1t_3}{2}-\frac{t_2^2}{4}+(-3t_1t_2+\frac{9t_3}{2})(-a) 
\end{array}\right), 
$$
\end{tiny}
If we express $\mathsf{R}_{\tau}$ as $u_a \frac{\partial}{\partial a}+u_b \frac{\partial}{\partial b}+u_c \frac{\partial}{\partial c}+u_{t_1} \frac{\partial}{\partial t_1}+u_{t_2} \frac{\partial}{\partial t_2}$, then the condition (\ref{R1=tau}) will lead to the equality
\[
\frac{1}{\Delta}(A_a u_a +A_b u_b +  A_c u_c + A_{t_1} u_{t_1} +A_{t_2} u_{t_2}) = \left( \begin{matrix}
0 & 0 & 0 \\
0 & 0 & -1 \\
0 & 0 & 0
\end{matrix}\right).
\]
Then after solving the equation, one may get $u_a, u_b, u_c, u_{t_1}, u_{t_2}.$ Similarly one can get $v_a, v_b, v_c, v_{t_1}, v_{t_2}$. $\hfill\square$
\section{Generalized period domain}
\label{gpd2021}
In this section we slightly modify the concept of generalized period domain introduced in \cite[Section 1]{ho18}, see also  \cite[Chapter 8, 11]{M2020}, to our context of mixed Hodge structures. This period domain comes with an action of a discrete group (the homology group) and the action of the algebraic group $\mathsf{G}$ described in Subsection \ref{alggrp}. The classical Griffiths period domain is obtained as the quotient by this algebraic group. There the values of the period map are defined, and it is a biholomorphism. Using the period map we finish the proof of Theorem \ref{23june2021}.

\subsection{Homology group}

Let $X$ be an abelian variety of dimension $g$ with two marked points $Y=\{O,P\}$. Writing $X=\C^g/\Lambda$, we can produce a natural basis $\delta_1,\ldots,\delta_{2g}\in H_1(X,\Z)$ induced by the generators of $\Lambda$. Using the wedge structure of $H_*(X,\Z)$, we can produce all the other homology groups. Since $H_1(X,\Z)\subseteq H_1(X,Y,\Z)$, we can add to this basis some $\delta_0\in H_1(X,Y,\Z)- H_1(X,\Z)$ in order to get a basis for all $H_1(X,Y,\Z)$ (it is a path connecting $O$ with $P$). For $X$ a principally polarized abelian variety we can take this basis such that its intersection matrix is given by
$$
\Psi=\begin{pmatrix}0 & 0 & 0 \\ 0 & 0 & -I_g \\ 0 & I_g & 0\end{pmatrix}.
$$
Denote $\Gamma_\Z:=\text{Aut}(H_1(X,Y,\Z),\langle\cdot,\cdot\rangle)$ and we have
$$
\Gamma_\Z=\text{Sp}(2g,\Z)\ltimes\Z^{2g}=\left\{\begin{pmatrix}1 & u \\ 0 & A\end{pmatrix}: \ A\in \text{Sp}(2g,\Z) \ , \ u\in \Z^{2g}\right\}.
$$
\subsection{Generalized period domain}
\label{sec8.2}
Now we introduce our generalized period domain $\mathsf{\Pi}$. Since we want it to be the domain where all periods live, it is natural to define it as the set of period matrices of the form
$$
\mathsf{P}:=\left[\int_{\delta_i}\alpha_j\right]=\begin{pmatrix}-1 & z_1 & z_2\\ 0 & x_1 & x_2 \\ 0 & x_3 & x_4\end{pmatrix},
$$
where $x_i\in\text{Mat}_{g\times g}$ and $z_i\in\C^g$ with
\begin{equation}
z_1\notin \Z^gx_1\oplus \Z^gx_3.        
\end{equation}
Note that the periods of $\alpha_0$ are zero over all closed paths. The Poincar\'e duality on $\mr{Gr}^W_1H^1_\dR((X,Y)/\C)\cong H^1_\dR(X/\C)$ translates into the equation
\begin{equation}
\label{16july2021}
\begin{pmatrix} 0 & -I_g \\ I_g & 0\end{pmatrix}=\begin{pmatrix} x_1^{\mathsf{tr}} & x_3^{\mathsf{tr}} \\ x_2^{\mathsf{tr}} & x_4^{\mathsf{tr}}\end{pmatrix}\begin{pmatrix} 0 & -I_g \\  I_g & 0\end{pmatrix}\begin{pmatrix} x_1 & x_2 \\  x_3 & x_4\end{pmatrix}=\begin{pmatrix} x_3^{\mathsf{tr}} & -x_1^{\mathsf{tr}} \\ x_4^{\mathsf{tr}} & -x_2^{\mathsf{tr}}\end{pmatrix}\begin{pmatrix} x_1 & x_2 \\ x_3 & x_4\end{pmatrix}
\end{equation}
$$
=\begin{pmatrix}
x_3^{\mathsf{tr}}x_1-x_1^{\mathsf{tr}}x_3 & x_3^{\mathsf{tr}}x_2-x_1^{\mathsf{tr}}x_4 \\ x_4^{\mathsf{tr}}x_1-x_2^{\mathsf{tr}}x_3 & x_4^{\mathsf{tr}}x_2-x_2^{\mathsf{tr}}x_4 
\end{pmatrix}.
$$
Hence the generalized period domain is given by the relations
\begin{equation}
x_1^{\mathsf{tr}}x_3=x_3^{\mathsf{tr}}x_1 \ , \ \ \ x_2^{\mathsf{tr}}x_4=x_4^{\mathsf{tr}}x_2 \ , \ \ \  x_1^{\mathsf{tr}}x_4-x_3^{\mathsf{tr}}x_2=I_g, 
\end{equation}
\begin{equation}
\label{Riemannrel}
\sqrt{-1}(x_3^{\mathsf{tr}}\overline{x_1}-x_1^{\mathsf{tr}}\overline{x_3}) \text{ is a positive matrix},     
\end{equation}
where \eqref{Riemannrel} corresponds to the second Hodge-Riemann bilinear relation. 
The corresponding Griffiths period domain $\mathsf{D}$ is obtained as the quotient under the action of the algebraic group $\mathsf{G}$. This is
\begin{equation}
\label{Griffpd}
\mathsf{D}:=\mathsf{\Pi}/\mathsf{G}\cong \C^g\times \uhp_g-\{(z,\tau) \ \big{|} \ \tau\in \uhp_g, \ z\in \Lambda_\tau=\Z^g\tau\oplus\Z^g\}
\end{equation}
given by
$$
\begin{pmatrix}
-1 & z_1 & z_2 \\ 0 & x_1 & x_2 \\ 0 & x_3 & x_4
\end{pmatrix}\mapsto (z_1x_3^{-1},x_1x_3^{-1}).
$$
Note that for $k'=-x_3^{-1}x_4x_3^\mathsf{tr}$, we have
$$
\begin{pmatrix}
-1 & z_1 & z_2 \\ 0 & x_1 & x_2 \\ 0 & x_3 & x_4
\end{pmatrix}\begin{pmatrix}
1 & 0 & -z_1k'-z_2x_3^{\mathsf{tr}} \\ 0 & x_3^{-1} & k' \\ 0 & 0 & x_3^{\mathsf{tr}}
\end{pmatrix}=\begin{pmatrix}
-1 & z_1x_3^{-1} & 0 \\ 0 & x_1x_3^{-1} & -I_g \\ 0 & I_g & 0
\end{pmatrix}.
$$
\subsection{The period map}
In the definition of the moduli space $\mathsf{T}$, the intersection matrix of $\alpha_1,\alpha_2,\ldots,\alpha_{2g}$ is the classical symplectic matrix.   This is minus the left hand side of \eqref{16july2021}. Therefore, we define the period map by inserting some minus sign: 
the period map is defined as
\[
\mathsf{P}: \mathsf{T} \to \mathsf{U}, \ \ 
\]
\[
t\mapsto 
(2\pi i)^{-\frac{1}{2g}}
\begin{bmatrix}
(2\pi i)^{\frac{1}{2g}}
\int_{\delta_j}\alpha_0 & \int_{\delta_j}\alpha_1 &\cdots &
\int_{\delta_j}\alpha_{g} & -\int_{\delta_j}\alpha_{g+1} &
\cdots &
-\int_{\delta_j}\alpha_{2g} 
\end{bmatrix},
\]
where $\mathsf{U}=\Gamma_\Z\backslash\mathsf{\Pi}$. The issue of the sign problem in \cite{M2012} is solved by using  $(-2\pi i)^{\frac{1}{2}}$. 

\begin{prop}
\label{isoperiodmap}
The period map $\mathsf{P}$ is a biholomorphism of complex manifolds.
\end{prop}

\begin{proof}
If we denote by $\mathsf{T}_0$ the moduli of enhanced principally polarized abelian varieties (without the marked points) and 
$$
\mathsf{P}_0:\mathsf{T}_0\rightarrow \mathsf{U}_0
$$
its corresponding period map. It is known that $\mathsf{P}_0$ is a biholomorphism, see for instance \cite{M2020, Fonseca}.  
Moreover, we have the following diagram
\begin{equation}
\begin{xy}
(100,0)*+{\mathsf{T} }="v1";(130,0)*+{\mathsf{U}}="v2";
(100,-20)*+{\mathsf{T}_0}="v3";(130,-20)*+{\mathsf{U}_0}="v4";
{\ar@{->}^{\mathsf{P}} "v1";"v2"};{\ar@{->}^{\mathsf{P}_0} "v3";"v4"};
{\ar@{->}^{\pi} "v1";"v3"};{\ar@{->}^{\rho} "v2";"v4"},
\end{xy}
\end{equation}
where $\pi$ and $\rho$ are the natural projection maps. Hence, to see that $\mathsf{P}$ is also a biholomorphism it is enough to show that it is a biholomorphism restricted to the fibers of $\pi$ and $\rho$. For every $t_0\in \mathsf{T}_0$
$$
\pi^{-1}(t_0)\cong (X_{t_0}-\{O\})\times\C^g
$$
On the other hand, if 
$$
\mathsf{P}_0(t_0)=\begin{pmatrix}
x_1 & x_2 \\ x_3 & x_4
\end{pmatrix}
$$
then we can identify
$$
\rho^{-1}(\mathsf{P}_0(t_0))=\Gamma_\Z\backslash\{(z_1,z_2)\in \C^g\times\C^g \ \big{|} \  z_1\notin \Z^gx_1\oplus\Z^gx_3 \}
$$
$$
\cong ((\C^g/(\Z^gx_1\oplus \Z^gx_3))-\{O\})\times\C^g.
$$
And the result follows noting that under the above identifications $\mathsf{P}=AJ\times id_{\C^g}$, where
$$
AJ:P\in X_{t_0}\mapsto z_1=\left(\int_{O}^P\alpha_1,\ldots,\int_O^P\alpha_g\right)\in \C^g/(\Z^gx_1\oplus \Z^gx_3)
$$
is just the usual Abel-Jacobi isomorphism.
\end{proof}

\subsection{$\tau$-map}
\label{23j2021}
In this subsection we will introduce the key map which explains our choice of the constant Gauss-Manin matrices of Theorem \ref{23june2021}. This is the $\tau$-map, which corresponds to a section of the quotient $\mathsf{\Pi}/\mathsf{G}$. In our context it is defined as 
$$
\tilde \tau: \mb{H}_g \times \mb{C}^g-\{(\tau,z) \ \big{|} \ \tau\in \uhp_g, \ z\in \Lambda_\tau=\Z^g\tau\oplus\Z^g\}\to \mathsf{\Pi},
$$
$$
(\tau,z)\mapsto  
\begin{bmatrix}
-1 & z & 0 \\
0 & \tau & -I_g \\
0 & I_g & 0
\end{bmatrix}.
$$
Since we have used the letter $\tau$ for $\tau\in\mb{H}_g$, we have named this map $\tilde\tau$. Its image is called the $\tau$-locus. 
Recalling the definition of the period map, it is elementary to see that
\begin{equation}
\label{GMpermatr}
\mr{d} \mathsf{P} = \mathsf{P} A^{\mr{tr}},   
\end{equation}
where $A$ is the Gauss-Manin matrix associated to the frame $[\alpha_0,\alpha_1,\ldots,\alpha_{2g}]$. Restricting \eqref{GMpermatr} to the $\tau$-locus we have 
$$
A =(\mathsf{P}^{-1} \mr{d}\mathsf{P})^{\mr{tr}} =  (\mr{d}\mathsf{P})^{\mr{tr}} \mathsf{P}^{-\mr{tr}}
$$
and 
\[
C_{ij}:= A\left(\frac{\partial}{\partial \tau_{ij}}\right) 
\]
which is the constant matrix such that all the entries are zero except $(i+1, g+j+1)$ and $(j+1, g+i+1)$ are $-1$.
Similarly, we have
\[
C_k :=  A\left(\frac{\partial}{\partial z_k}\right) 
\]
which is the constant matrix such that all the entries are zero except $(k+1, 1)$ is $-1$. In consequence, the vector fields determined by Theorem \ref{23june2021} are tangent to the $\tau$-locus. In other words, in the case of elliptic curves, the coordinates of $\mathsf{T}$ restricted to the $\tau$-locus, are solutions to the modular differential equations defined by $R_z$ and $R_\tau$ of Theorem \ref{23.06.2021}.
\subsection{Elliptic curves}
In this case the $\tau$-locus  is given by 
$
\left( \begin{matrix}
-1 & z & 0 \\
0 & \tau & -1 \\
0 & 1 & 0
\end{matrix} \right)$, where $ \tau\in\mf{h}, \ z\in\C-\Lambda_\tau$  and 
\[
A\left(\frac{\partial}{\partial \tau}\right) =
\left( \begin{matrix}
0 & 0 & 0 \\
0 & 1 & 0 \\
0 & 0 & 0
\end{matrix} \right)^{\mr{tr}}  \left( \begin{matrix}
-1 & 0 & z \\
0 & 0& 1 \\
0 & -1 & \tau
\end{matrix} \right)^{\mr{tr}} =\left( \begin{matrix}
0 & 0 & 0 \\
0 & 0 & -1 \\
0 & 0 & 0
\end{matrix} \right).
\]
Similarly, we have
\[
A\left(\frac{\partial}{\partial z}\right) =
\left( \begin{matrix}
0 & 1 & 0 \\
0 & 0 & 0 \\
0 & 0 & 0
\end{matrix} \right)^{\mr{tr}}  \left( \begin{matrix}
-1 & 0 & z \\
0 & 0& 1 \\
0 & -1 & \tau
\end{matrix} \right)^{\mr{tr}} =\left( \begin{matrix}
0 & 0 & 0 \\
-1 & 0 & 0 \\
0 & 0 & 0
\end{matrix} \right).
\]
\subsection{Proof of Theorem \ref{23june2021}}
The main idea behind is the same as the proof \cite[Theorem 11.5]{M2020}. We first observe that by Proposition \ref{isoperiodmap} the period map
\[
\mathsf{P}: \mathsf{T} \to \mathsf{U}
\]
 is a biholomorphism. It is enough to prove the existence and uniqueness of $\mathsf{v}_{ij}$ and $\mathsf{v}_k$ in the period domain $\mathsf{U}$. Since $A$ is defined over $\mb{Q}$, by uniqueness, it follows that these vector fields are also defined over $\mb{Q}$.
Note that $\mathsf{U} = \Gamma_{\mb{Z}} \backslash \mathsf{\Pi}$. Combining  \cite[Proposition 8.10]{M2020} and the definition of $C_{ij}, C_k$ in the previous subsection, 
it is enough to check the equality
\[
C_{ij}^{\mr{tr}} \Psi + \Psi C_{ij} = 0 \ , \ \ \  C_{k}^{\mr{tr}} \Psi + \Psi C_{k} = 0
\]
which can be proved by the direct computation. For the last statement on Lie brackets of these vector fields, it follows from \cite[Proposition 6.17]{M2020} 
and the fact that the Lie bracket of matrices $C_{ij}, C_k$'s is zero.\qed
\section{Quasi Jacobi forms of index zero}
\label{sect9}
In this last section we use the period map $\mathsf{P}$ and the $\tau$-map described in the previous section, to construct the map $\tmap:\mathsf{D}\rightarrow\mathsf{T}$ announced in the introduction. Using this map we lift the solutions of the modular differential equations on $\mathsf{T}$ to $\mb{H}_g \times \mb{C}^g-\{(\tau,z) \ \big{|} \ \tau\in \uhp_g, \ z\in \Lambda_\tau=\Z^g\tau\oplus\Z^g\}$. Then, using both actions of the discrete group $\Gamma_\Z$ and the action of the algebraic group $\mathsf{G}$, we get the modularity conditions on each solution of the modular vector fields, revealing to us which quasi Jacobi form of index zero corresponds to each solution.

\subsection{Quasi Jacobi forms of index zero}
\label{sec9.1}
Recall that there is a canonical morphism of schemes
\[
\mathsf{T} \to \mathsf{S}, \ \ \ (X, Y, [\alpha_0, \alpha_1, \ldots, \alpha_{2g}]) \to 
(X, Y, \alpha_1 \wedge \alpha_2 \wedge \cdots \wedge \alpha_{g}).
\]
The map is surjective and hence we have an inclusion $k[\mathsf{S}] \subset k[\mathsf{T}]$.
\begin{df}
The algebra of algebraic Jacobi (resp. quasi Jacobi) forms of index zero is by definition $k[\mathsf{T}]$ (resp. $k[\mathsf{S}]$). A function $f \in k[\mathsf{S}]$ is called an algebraic Jacobi-form of weight $k$ and index zero if 
\[
f(t \bullet g) = f(t) g^{-k}, \forall g \in \mb{G}_m.
\]
\end{df}
Over the field of complex numbers, we consider the composition $\tmap$ of the maps:
\begin{equation}
\label{14ju2021}
\tmap: \mb{H}_g \times \mb{C}^g -\{(\tau,z) \ \big{|} \ \tau\in \uhp_g, \ z\in \Lambda_\tau=\Z^g\tau\oplus\Z^g\}\to \mathsf{\Pi} \to \mathsf{U} \xrightarrow{\mathsf{P}^{-1}} \mathsf{T},
\end{equation}
where the first map is the $\tau$-map and the second map is the canonical quotient map. For the definition of this map without assuming that the period map is a biholomorphism see 
\cite[Section 8.5]{M2020}. 
\begin{df}
The algebra of holomorphic Jacobi (resp. quasi Jacobi) forms of index zero is by definition the pull-back of the algebra of algebraic Jacobi (resp. quasi Jacobi) forms of index zero by the map $\mb{H}_g \times \mb{C}^g\to \mathsf{T}$.
\end{df}
We can show that a holomorphic Jacobi form $\phi$ of weight $k$ and index zero satisfy the functional equations:
\begin{enumerate}
\item
$\phi(M\tau, z(C\tau +D)^{-1}) = \mr{det}(C\tau+D)^k 
\phi(\tau, z)$, where
$M=\begin{pmatrix}A&B\\ C&D\end{pmatrix} \in Sp(2g, \mb{Z})$,
\item
$\phi(\tau, z+\lambda\tau+\mu) = 
\phi(\tau, z)$, where $(\lambda, \mu) \in \mb{Z}^{2g}$.
\end{enumerate}
The second functional equation implies that holomorphic Jacobi forms which are also holomorphic in the lattice points 
$\Lambda_\tau$ are necessarily constant in $z$, and so they are classical quasi modular forms.  

\subsection{Classical generators}
We recall some classical quasi Jacobi forms of index zero. The first typical cases are classical modular forms, for example, the Eisenstein series 
\begin{equation} \label{E_i}   
E_{2i} = 1+b_i \sum^{\infty}_{n=1}\left(\sum_{d|n} d^{2i-1}\right)q^n \ , \ \ \ \  (b_1,b_2,b_3) = (-24,240,-504),
\end{equation}
where $q = e^{2\pi i \tau}$.
Let $(\tau, z) \in \mb{H} \times \mb{C}$ and let $y= - e^{2\pi i z}$ and $q= e^{2\pi i \tau}$. We let
\begin{equation} \label{F}
F(\tau, z) = \frac{\theta_1(\tau, z)}{\eta^3(\tau)} = (y^{1/2}+y^{-1/2})\prod_{m \geq 1} \frac{(1+yq^m)(1+y^{-1}q^m)}{(1-q^m)^2}
\end{equation}
where the expansion is considered in the region $|y|<1$, and 
the logarithmic derivative
\begin{equation} \label{J}
J_1(\tau, z) =  y\frac{\mr{d}}{\mr{d}y} \log F(y, q)
\end{equation}
and the Weierstrass elliptic function $\wp(\tau, z)$ together with its derivative
\begin{equation} \label{p'}
\wp^{'}(\tau, z) = y\frac{\mr{d}}{\mr{d}y} \wp(\tau, z).
\end{equation}
see \cite[Appendix B]{Ober2018}.
We remark that $E_{2k}(\tau), \wp(\tau, z), \wp^{'}(\tau, z)$
are index zero Jacobi forms and  $J_1(\tau, z)$ is a quasi Jacobi form. 
Moreover they can have a pole in the fundamental region $\{x+y\tau \mid 0 \leq x, y < 1\}$
only at $z=0$.  

\subsection{Proof of Theorem \ref{2021china}}
The statement for $t_1,t_2,t_3$ follows from \cite[Section 6.4]{M2012}.
We compute the transformation law for $\tilde a:=a\circ \tmap, \tilde b:=b\circ \tmap, \tilde c:=c\circ \tmap, \tilde t_i:=t_1\circ\tmap, \ \ i=1,2,3$ at first. 
Note that the difference between the period map in this text and the mentioned reference is given by the action of the element $\mat{i^{-1}}{0}{0}{i}$ of the algebraic group used in this reference. This transforms $(t_1,t_2,t_3)$ into $(-t_1,t_2,-t_3)$. 
Let $A = \left(\begin{matrix}
1 & \lambda & \mu \\
0 & \alpha & \beta \\
0 & \gamma & \delta
\end{matrix}\right) \in SL_2(\mb{Z}) \ltimes \mb{Z}^2$.
Then we have an equality
\begin{equation} \label{key equality}
\begin{split}
&\left(\begin{matrix}
1 & \lambda & \mu \\
0 & \alpha & \beta \\
0 & \gamma & \delta
\end{matrix}\right)
\left(\begin{matrix}
-1 & z & 0 \\
0 & \tau & -1 \\
0 & 1 & 0
\end{matrix}\right) \\
=& \left(\begin{matrix}
-1& \frac{z+\lambda\tau+\mu}{\gamma\tau+\delta} & 0\\
0 & \frac{\alpha\tau+\beta}{\gamma\tau+\delta} & -1 \\
0 & 1 & 0
\end{matrix}\right)
\left(\begin{matrix}
1  & 0 & R_A(\tau, z) \\
0 & \gamma\tau+\delta & -\gamma \\
0 & 0 & \frac{1}{\gamma\tau+\delta}
\end{matrix}\right),
\end{split}
\end{equation}
where $g = \left(\begin{matrix}
1 & 0 & R_A(\tau, z) \\
0 & \gamma\tau+\delta & -\gamma \\
0 & 0 & \frac{1}{\gamma\tau+\delta}
\end{matrix}\right) \in \mathsf{G}$ and
$R_A(\tau, z) =  \frac{\lambda \delta - \mu \gamma - z \gamma}{\gamma \tau+\delta}.$ 
We have 
\begin{equation} \label{transformation law}
\begin{split}
&\widetilde{a}(\tau, z) = a(i(\tau, z)) =  a\left(i\left(\frac{\alpha\tau+\beta}{\gamma\tau+\delta}, \frac{z+\lambda \tau+\mu}{\gamma\tau+\delta}\right)g\right) \\
= & (a \bullet g)\left(i\left(\frac{\alpha\tau+\beta}{\gamma\tau+\delta}, \frac{z+\lambda \tau+\mu}{\gamma\tau+\delta}\right)\right) = (\gamma\tau+\delta)^{-2} \widetilde{a}\left(\frac{\alpha\tau+\beta}{\gamma\tau+\delta}, \frac{z+\lambda \tau+\mu}{\gamma\tau+\delta}\right)
\end{split} 
\end{equation}
Hence $\widetilde{a}$ satisfies the first two properties of weak Jacobi forms for $k=2, m = 0$. 
For the transfomation law of $\tilde{c}$, consider the following three type generators of ${\rm SL}(2,\mb{Z}) \ltimes \mb{Z}^2$:
\[
H= \left(\begin{matrix}
1 & \lambda & \mu \\
0 & 1 & 0 \\
0 & 0 & 1
\end{matrix}\right), \  T =\left(\begin{matrix}
1 & 0 & 0 \\
0 & 1 & 1 \\
0 & 0 & 1
\end{matrix}\right), \ S = \left(\begin{matrix}
1 & 0 & 0 \\
0 & 0 & -1 \\
0 & 1 & 0 
\end{matrix}\right).
\]
Using the equality (\ref{key equality}), one may find that
\[
R_H = \lambda, \ \ R_T = 0, \ \ R_S = -\frac{z}{\tau}.
\]
Recall that $g \in \mathsf{G}$ acts on $c$ by $c \bullet g = k^{-1}c + v$.
Similar considerations (equality (\ref{transformation law})) as above yields that
\begin{eqnarray}
\widetilde{c}(\tau, z+\lambda \tau + \mu) &= & \widetilde{c}(\tau, z) - \lambda \\
\widetilde{c}(\tau+1, z) & =& \widetilde{c}(\tau, z) 
\\
\widetilde{c}\left(-\frac{1}{\tau}, \frac{z}{\tau}\right) & = & \tau \widetilde{c}(\tau, z) + z
\end{eqnarray}
In conclusion, the transformation law for $a,b, c,t_1,t_2,t_3$ coincide with the transformation law for quasi Jacobi forms $\wp$ (resp. $\wp^{'}$, $J_1$, $E_2, E_4, E_6$) of weight $2$ (resp. $3, 1, 2, 4, 6$) and index zero. Note that $\tilde a$ and $\tilde b$ are holomorphic except at the points $(z, \tau)$, where $z \in \Lambda_{\tau}$. In other words, for any fixed $\tau$, $z= 0$ is a pole of $\tilde a$ and $\tilde b$. Because $R_z(a) = b$ and $b^2 = 4a^3 -t_2a -t_3$, the functions $\tilde a$ and $\tilde b$ are elliptic functions of pole order respectively $2$ and $3$ at $z=0$, and hence, up to constant they are equal to $\wp$ and $\wp^{'}$. The computation of those constants also follows from $b^2=4a^3-t_2a-t_3$. It remains to prove the theorem for $\tilde c$.  
Through (\ref{R1=tau}) and (\ref{R2=z}), we know that 
$\tilde a,\tilde b, \tilde c,\tilde t_1,\tilde t_2,\tilde t_3$ form a solution of $-R_\tau$ and $R_z$.     
These are the same differential equations satisfied by $J_1, \wp$ and $\wp^{'}, E_2,E_4$ computed in \cite[Lemma 48]{Ober2018} (after inserting the $2\pi i$ factors). 
It is easy to see that if we have two solutions of $-R_\tau$ and $R_z$ with the same $t_1,t_2,a,b$ then their $c$ is also the same. $\hfill\square$

\subsection{The Serre-Jacobi derivative}
The algebra of classical modular froms is not closed under the derivative $\frac{\partial}{\partial \tau}$, but the modularity can be recovered by adding a multiple of $E_2$. We then get a differential operator on modular forms, which is called the Serre derivative. If we identify the algebra of quasi modular forms with $\mb{C}[t_1, t_2, t_3]$, the Serre derivative can be defined as
\[
\partial^S: f \to \partial_{\tau} f - (w_f-s_f) t_1 f, 
\]
where $w_f$ is the weight of $f$ and $s_f$ is the depth or differential order (degree in $t_1$) of $f$ . For example, Ramanujan differential equations, see \cite[Proposition 4.1]{M2012}, is equivalent to
\[
\partial^S(t_1) =- \frac{t_2}{12}, \ \ \partial^S(t_2) = -6t_3, \ \ \partial^S(t_3)=-\frac{t_2^2}{3}.
\]
For even weight Jacobi forms, Oberdieck  \cite[Lemma 10]{Ober2012} defines the similar differential operator, which is called the Jacobi-Serre derivative.  We slightly modify such a diffential operator as follows. 
We first note that a solution of $R_\tau$ and $R_z$ is given by 
$$
-(2\pi i)\wp(\tau, z) \ , \ \  i^{\frac{3}{2}}(2\pi i)^{\frac{3}{2}}\wp^{'}(\tau, z) \ , \ \ - i(2\pi i)^{\frac{1}{2}}J_1(\tau, z),
$$
$$
\frac{2\pi i}{12}E_2(\tau) \ , \ \ 12\left(\frac{2\pi i}{12}\right)^2E_4(\tau) \ , \ \ 8\left(\frac{2\pi i}{12}\right)^3E_6(\tau),
$$
which we denote by $a,b,c,t_1,t_2,t_3$. We define 
\[
\partial^J: \mb{C}[a,b,c,t_1,t_2,\frac{1}{\Delta}] \to \mb{C}[a,b,c,t_1,t_2,\frac{1}{\Delta}]
\]
as follows:
\[
\partial^J(f) = \partial_{\tau}f - (w_f-s_f) t_1 f - c \partial_zf,
\]
for $f \in \mb{C}[a,b,c,t_1,t_2,\frac{1}{\Delta}]$ , where $k_f$ is the weight of $f$ and $s_f$ is the degree in $t_1$ variable. 
\begin{prop}
 The Serre-Jacobi derivative $\partial^J(f)$ does not increse the degree in $t_1$ and $c$. It increases the weight by $2$. 
\end{prop}
\begin{proof}
First, note that
$$\partial^J(a)= -2a^2 + \frac{t_2}{3},  \ \ \partial^J(b) =   -3ab,  \ \ \partial^J(c) =  - \frac{b}{2}-ct_1,\ \ \partial^J(t_1)=-\frac{t_2}{12},
\partial^J(t_2) = -6t_3
$$
We need to consider polynomials of the form $f:=t_1^nP(a,b,c,t_2)$ (resp. $f:=c^nP(a,b,t_1,t_2)$), where $P$ is homogeneous,   
and observe that $\partial^J(f)$ does not increase degree in $t_1$ (resp. $c$). The computation is easy and it is left to the reader. 
\end{proof}

\bibliographystyle{alpha}
\bibliography{biblio}

\address{Jin Cao\\
Yau Mathematical Sciences Center, 
Tsinghua University, Beijing, China}
\

\email{\texttt{caojin@mail.tsinghua.edu.cn}}

\bigskip

\address{Hossein Movasati\\
Instituto de Matem\'atica Pura e Aplicada, IMPA, Estrada Dona Castorina, 110, 22460-320, Rio de Janeiro, RJ, Brazil}
\

\email{\texttt{hossein@impa.br}}

\bigskip

\address{Roberto Villaflor Loyola\\
Facultad de Matemáticas, Pontificia Universidad Católica de Chile, Campus San Joaquín, Avenida Vicuña Mackenna 4860, Santiago, Chile}
\

\email{\texttt{roberto.villaflor@mat.uc.cl}}

\end{document}